\documentclass[]{article}
\usepackage{geometry}
\usepackage{amsmath,amssymb,amsfonts,amsthm,amsopn,dsfont, mathrsfs}
\usepackage{bm}
\usepackage{graphicx}
\usepackage{float}

\usepackage{standalone}
\usepackage{tikz}
\usepackage{pgfplots}
\usepackage{pgfplotstable}
\pgfplotsset{compat=newest}
\usepgfplotslibrary{patchplots}
\usepackage{subcaption}
\usepackage{colortbl}
\usepackage{booktabs}

\usepackage[T1]{fontenc}

\usepackage{xcolor}

\usepackage[ruled,vlined,linesnumbered]{algorithm2e}
\usepackage{multirow}

\newcommand{\hsp}{\mathcal{H}^{\Gamma}}

\newcommand{\extg}{\mathcal{E}_{_{\Gamma}}}
\newcommand{\exts}{\mathcal{E}_{_{\Sigma}}}
\newcommand{\trg}{\gamma_{_{\Gamma}}}

\title{Extended Finite Elements for 3D-1D coupled problems via a PDE-constrained optimization approach}

\author{Denise Grappein\footnote{denise.grappein@polito.it} , 
	Stefano Scial\'o\footnote{stefano.scialo@polito.it} , 
	Fabio Vicini\footnote{fabio.vicini@polito.it}\medskip\\{\footnotesize Dipartimento di Scienze Matematiche, Politecnico di Torino,
		Corso Duca degli Abruzzi 24,}\\{\footnotesize 10129 Torino, Italy. Members of INdAM research group GNCS.}
}

\begin{document}
\maketitle

\begin{abstract}
In this work, we propose the application of the eXtended Finite Element Method (XFEM) in the context of the coupling between three-dimensional and one-dimensional elliptic problems. In particular, we consider the case in which the 3D-1D coupled problem arises from the geometrical model reduction of a fully three dimensional problem, characterized by thin tubular inclusions embedded in a much wider domain. 
In the 3D-1D coupling framework, the use of non conforming meshes is widely adopted. However, since the inclusions typically behave as singular sinks or sources for the 3D problem, mesh adaptation near the embedded 1D domains may be necessary to enhance solution accuracy and recover optimal convergence rates. An alternative to mesh adaptation is represented by the XFEM, which we here propose to enhance the approximation capabilities of an optimization-based 3D-1D coupling approach. An effective quadrature strategy is devised in order to integrate the enrichment functions and numerical tests on single and on intersecting segments are proposed to demonstrate the effectiveness of the approach.
\end{abstract}

\section{Introduction}\label{sec:intro}
Coupled partial differential equation problems on 3D and 1D domains arise from the application of dimensional reduction models to equi-dimensional problems where cylindrical or nearly-cylindrical inclusions with small cross sectional size are embedded in a larger 3D domain, \cite{Koppl2018,Zunino2019,BGS3D1D2022}. 
The treatment of such narrow and elongated regions as one-dimensional manifolds reduces the overhead in simulations related to the generation of a computational mesh inside the inclusions. Suitable matching conditions need to be added at the interfaces to close the problem, depending on the nature of the described physical phenomenon: in some cases the solution is expected to be continuous at the interface, as in the description of damaged vessels in tumour induced angiogenesis \cite{Giverso2016}, of thin membranes \cite{Chaplain2019}, or in the modeling of fiber reinforced materials \cite{Steinbrecher2020}; in other cases, filtration like conditions, yielding a discontinuity at the 3D-1D interfaces are preferred, as for plant-roots nutrient uptake from the subsoil \cite{Schroder2012,Koch2018}, in geological applications \cite{Gjerde2018a}, or again in angiogenesis \cite{Notaro2016,Koppl2020}.
However, the mathematical formulation of 3D-1D coupled problems requires non-standard approaches, and specialized numerical schemes are needed to correctly account for the presence of singularities. A possibility lies in the use of regularizing functions to approximate the singular terms \cite{Heltai2019,Koch2020}, or of lifting techniques \cite{Koppl2018}. In \cite{Gjerde2019} the solution is split in a regular part, approximated by standard methods, and an irregular part, for which an analytical solution is given. Domain decomposition approaches are finally proposed in \cite{Kuchta2021} based on Lagrange multipliers, and in \cite{BGS3D1D2022,BGS3D1Ddisc} where a PDE-constrained optimization method is presented. 

The use of a 3D mesh non conforming to the 1D domains is quite standard.
However, in some cases, sub-optimal convergence rates are observed unless adaptive refinement close to the singularity is used, see e.g. \cite{Koppl2020,Heltai2019}. 
In this work we adopt the eXtended Finite Element Method (XFEM) as an alternative to mesh refinement. The application of XFEM to 3D problems with singular sources has been proposed in \cite{Gracie2010}, in particular for quasi 3D problems describing the effect of well leakage in aquifers. In \cite{BREZINA2021} the methodology was extended to fully 3D-1D coupled problems in mixed formulation. The function space for the velocity variable is enriched, and non intersecting segments entirely crossing the computational domain are considered in the numerical examples.

Here, we focus on the application of the XFEM to enhance the approximation capabilities of the optimization-based 3D-1D coupling strategy proposed in \cite{BGS3D1D2022,BGS3D1Ddisc}. Such method is based on a three-field domain decomposition strategy, in which additional interface variables are introduced to de-couple the problem on the inclusions from the problem in the bulk domain. A cost functional is introduced to measure the error in satisfying the desired matching condition at the interfaces, and minimized to recover a global solution. Different interface conditions are considered in \cite{BGS3D1D2022} and in \cite{BGS3D1Ddisc}, resulting in two different formulations of the method. In the present work we consider flux conservation and pressure continuity at the interface as in \cite{BGS3D1D2022}. However, the proposed approach can be easily extended to other interface conditions, such as the ones considered in \cite{BGS3D1Ddisc}, or even to different formulations of the problem. We enrich the function space of the 3D pressure variable with a globally continuous function, having a log-like behavior outside the inclusion and being instead constant inside it. The choice of the enrichment function is based on the results provided in \cite{Gjerde2019}, adapted to the present case. We suggest an ad-hoc quadrature scheme for the numerical integration of the resulting irregular basis functions. We consider intersecting/branching inclusions, possibly ending inside the domain. 

The manuscript is organized as follows: the model problem is presented in Section \ref{not_modprob}, and an overview on the PDE-constrained optimization approach is provided in Section \ref{sec:opt_form}; Section \ref{discr} is devoted to the general discretization of the optimization problem, while the details on the application of the XFEM are provided in Section \ref{xfem_func}. In Section \ref{quad} we propose a quadrature strategy, suitably designed to integrate the enriched basis functions and finally, in Section \ref{num_res}, some numerical experiments are presented, in order to validate the proposed approach.

\section{Notation and model problem}\label{not_modprob}
\begin{figure}
	\centering
	\includegraphics[width=0.37\textwidth]{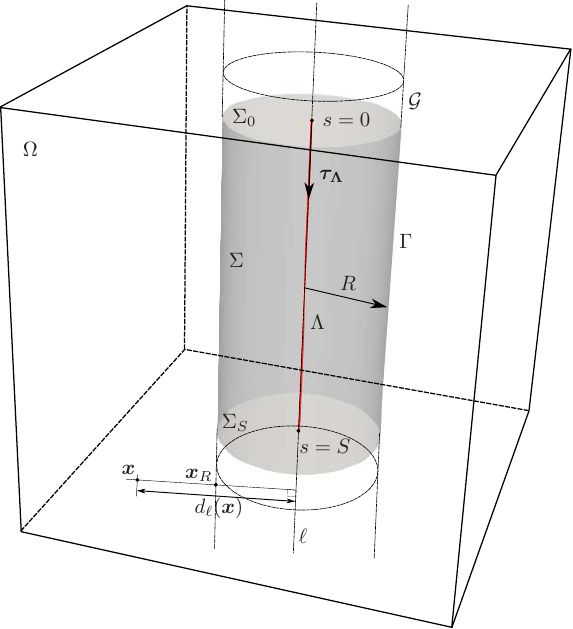}
	\caption{Domain with single inclusion and description of notation. The size of the inclusion is exaggerated for description purposes.}
	\label{fig:notation}
\end{figure}

We consider a convex domain $\Omega \subset \mathbb{R}^3$, characterized by the presence of a thin cylindrical inclusion $\Sigma \subset \Omega$ of radius $R$, see Figure~\ref{fig:notation}. We assume $R$ to be smaller than both the diameter of the domain and the length of the cylinder. The centerline of $\Sigma$ is denoted by $\Lambda=\left\lbrace \bm{\lambda}(s),~s \in (0,S)\right\rbrace $, and $\bm{\tau_\Lambda}$ is the unit tangent vector to $\Lambda$. We further call $\ell$ the line passing through the centreline of $\Sigma$ and $\mathcal{G}$ the lateral surface of an infinite cylinder with centreline $\ell$ and radius $R$.
The boundary of $\Omega$ is denoted by $\partial \Omega$ and is split into two subsets: the Dirichlet boundary $\partial \Omega_{\mathrm{d}}$ and the Neumann boundary $\partial \Omega_{\mathrm{n}}$, such that $\partial \Omega=\overline{\partial \Omega_{\mathrm{d}}}\cup \overline{\partial \Omega_{\mathrm{n}}}$, with $\partial \Omega_{\mathrm{d}}\cap \partial \Omega_{\mathrm{n}}=\emptyset$ and $|\partial \Omega_{\mathrm{d}}| > 0$. The boundary of $\Sigma$ is instead split into the lateral surface $\Gamma$ and the two end sections $\Sigma_0$ and $\Sigma_S$, i.e. $\partial \Sigma=\overline{\Gamma}\cup \overline{\Sigma_0}\cup \overline{\Sigma_S}$. The symbol $\Sigma(s)$, $s \in (0,S)$, is used to denote a generic cross-section of $\Sigma$.
Finally, we define $D:=\Omega\setminus \overline{\Sigma}$, the domain without the inclusion, having boundary $\partial D=\partial \Omega \cup \partial \Sigma$.

Let us now consider a diffusion problem in $D$ and $\Sigma$ with unknown pressures $u \in D$ and $\tilde{u} \in \Sigma$:
\begin{minipage}[t]{0.5\textwidth}\vspace{-0.3cm}
	\begin{align}
	-&\nabla \cdot (K \nabla u)=f & \text{in } D\label{1}\\
	&u_{|_{\Gamma}}=\psi &~\text{ on } \Gamma\\
	&K\nabla u \cdot \bm{n}=\phi &~\text{on } \Gamma \\[0.5em]
	&u=0 &\text{on } \partial \Omega_{\mathrm{d}}\\
	&\nabla u \cdot \bm{n}=0 &\text{on } \partial \Omega_{\mathrm{n}}\\
	&\nabla u\cdot\bm{n}=0 &\text{ on } \Sigma_0 \cup \Sigma_S
	\end{align}		
\end{minipage}
\begin{minipage}[t]{0.5\textwidth}\vspace{-0.3cm}
	\begin{align}
	-&\nabla \cdot(\tilde{K}\nabla\tilde{u})=g&\text{ in } \Sigma\\
	&\tilde{u}_{|_{\Gamma}}=\psi &\text{ on } \Gamma\\
	&\tilde{K}\nabla \tilde{u}\cdot \tilde{\bm{n}}=-\phi &\text{ on } \Gamma\\[0.5em]
	&\nabla \tilde{u}\cdot\bm{\tilde{n}}=0 &\text{ on } \Sigma_0 \cup \Sigma_S \label{end}
	\end{align}\medskip\\
\end{minipage}
For the sake of simplicity, we assume that all boundary conditions, except the interface conditions prescribed on $\Gamma$, are homogeneous. 
The vector $\bm{n}$ denotes the outward-pointing unit normal to $\partial D$, while $\bm{\tilde{n}}=-\bm{n}$ is the outward pointing unit normal to $\partial \Sigma$.

\section{Optimization formulation for the 3D-1D reduced problem}\label{sec:opt_form}

Following \cite{BGS3D1D2022}, to which we refer for details, the above 3D-3D coupled problem is reformulated as a 3D-1D coupled problem through a suitable choice of function spaces for the solution. Given the small radius of the inclusion, the solution is assumed to be constant on its cross sections and their boundaries. Let us hence define two extension operators:
\begin{equation*}
\mathcal{E}_{_\Sigma}: H^1(\Lambda) \rightarrow H^1(\Sigma) ~\text{  and  }~\mathcal{E}_{_\Gamma}: H^1(\Lambda) \rightarrow H^{\frac{1}{2}}(\Gamma),
\end{equation*}
which, given a function $\hat{v} \in H^1(\Lambda)$, uniformly extend the value $\hat{v}(s)$, $s \in [0,S]$ to the cross section $\Sigma(s)$ of the cylinder, i.e. $\exts\hat{v}(\bm{x})=\hat{v}(s)~\forall \bm{x}\in \Sigma(s)$, and to the boundary $\Gamma(s)$ of $\Sigma(s)$, i.e. $\extg\hat{v}(\bm{x})=\hat{v}(s)~\forall \bm{x} \in \Gamma(s)$. Given the trace operator $\gamma_{_\Gamma}:H^1(D)\cup H^1(\Sigma)\rightarrow H^{\frac{1}{2}}(\Gamma)$, defined as $\gamma_{_\Gamma}v=v_{|_\Gamma} ~\forall v \in H^1(D)\cup H^1(\Sigma)$, and setting $\hat{V}=H^1(\Lambda)$, the following spaces are introduced:
\begin{displaymath}
\widetilde{V}=\lbrace v \in H_0^1(\Sigma): v =\exts\hat{v}, ~\hat{v} \in \hat{V} \rbrace, \quad \mathcal{H}^{\Gamma}=\lbrace v \in H^{\frac{1}{2}}(\Gamma): v =\extg\hat{v}, ~\hat{v} \in \hat{V} \rbrace
\end{displaymath}
$$V_D=\left\lbrace v \in H^1(D): v_{|_{\partial \Omega_{\mathrm{d}}}}=0 \text{ and } \gamma_{_\Gamma}v \in \mathcal{H}^{\Gamma}\right\rbrace,$$
whose functions satisfy the hypothesis on the regularity of the solution. 
Denoting by $(\cdot,\cdot)_\star$ the $L^2$-scalar product on a generic domain $\star$, by $X'$ the dual of a space $X$, and by $\langle\cdot,\cdot\rangle_{X', X}$ the duality pairing between the two spaces, the weak formulation of \eqref{1}-\eqref{end} reads: {\em find} $(u, \tilde{u}) \in V_D\times \widetilde{V}$, $\psi \in \mathcal{H}^{\Gamma} $, $\phi \in {\mathcal{H}^{\Gamma}}'$ \textit{such that}
\begin{align}
&(K\nabla u, \nabla v)_{D}-\left\langle \phi,\gamma_{_\Gamma}v \right\rangle_{{\mathcal{H}^{\Gamma}}', {\mathcal{H}^{\Gamma}}}=(f,v)_{D}~&\forall v \in V_D,~ \label{eq_var_u} \\
&(\tilde{K}\nabla \tilde{u}, \nabla \tilde{v})_{\Sigma}+\left\langle \phi,\gamma_{_\Gamma}\tilde{v} \right\rangle_{{\mathcal{H}^{\Gamma}}', {\mathcal{H}^{\Gamma}}}=(g,\tilde{v})_{\Sigma} &\forall \tilde{v} \in \widetilde{V}, \label{eq_var_uhat}\\
&\left\langle \gamma_{_\Gamma}u-\psi,\eta\right\rangle_{\mathcal{H}^{\Gamma},{\mathcal{H}^{\Gamma}}'}=0&~\forall \eta \in  {\mathcal{H}^{\Gamma}}', \label{condpsi_u}\\
&\left\langle \gamma_{_\Gamma}\tilde{u}-\psi,\eta\right\rangle_{\mathcal{H}^{\Gamma},{\mathcal{H}^{\Gamma}}'}=0&~\forall \eta \in  {\mathcal{H}^{\Gamma}}'\label{condpsi_hat}.
\end{align}
A well posed 3D-1D formulation follows by operating a geometrical reduction of the operators:
\begin{align*}
\left\langle \phi,\gamma_{_\Gamma}v \right\rangle_{{\mathcal{H}^{\Gamma}}', {\mathcal{H}^{\Gamma}}}&=\int_{\Gamma}\phi~\gamma_{_\Gamma}v~d\Gamma=\int_0^S\Big( \int_{\Gamma(s)}\phi~\gamma_{_\Gamma}v~dl\Big) ds=\\ &=\int_0^S|\Gamma(s)|\overline{\phi}(s)\check{v}(s)~ds=\left\langle |\Gamma|\overline{\phi},\check{v}\right\rangle_{\hat{V}',\hat{V}} \quad \forall v \in V_D,
\end{align*}
\begin{equation*}
(\tilde{K}\nabla \tilde{u}, \nabla \tilde{v})_{\Sigma}=\int_{\Sigma}\tilde{K}\nabla \tilde{u}\nabla\tilde{v}~d\sigma=\int_0^S\tilde{K}|\Sigma(s)|\cfrac{d\hat{u}}{ds}~\cfrac{d\hat{v}}{ds}~ds,
\end{equation*}
where $\check{v} \in \hat{V}$ is such that  $\gamma_{_\Gamma}v=\extg\check{v}$, and $\hat{u},\hat{v} \in \hat{V}$ such that $\tilde{u}=\exts\hat{u}$, $\tilde{v}=\exts\hat{v}$.
The quantities $|\Gamma(s)|$ and $|\Sigma(s)|$ are the measure of $\Gamma(s)$ and $\Sigma(s)$, respectively.

Instead of solving the coupled system of equations \eqref{eq_var_u}-\eqref{condpsi_hat}, we re-write it as a PDE-constrained optimization problem. This is done by introducing a cost functional to measure the error in fulfilling the coupling conditions \eqref{condpsi_u}-\eqref{condpsi_hat}, and looking at the solution as the minimum of this functional, constrained by the constitutive equations on the 3D and 1D domains:
\begin{eqnarray}
&&\min_{\overline{\phi},\hat{ \psi}} J=\cfrac{1}{2}\left( ||\gamma_{_\Gamma}u(\overline{\phi},\hat{ \psi})-\extg \hat{\psi}||_{\hsp}^2+||\gamma_{_\Gamma}\exts \hat{u}(\overline{\phi},\hat{ \psi})-\extg \hat{\psi}||_{\hsp}^2\right),	\label{minJ} \\
&&\text{such that, for } \check{v}, \hat{\psi} \in \hat{V}: \trg v=\extg\check{v} \ \text{and} \ \psi=\extg\hat{\psi}: \nonumber \\
&&(K\nabla u, \nabla v)_{D}+\alpha(|\Gamma|\check{u},\check{v} )_{\Lambda}-\left\langle |\Gamma|\overline{\phi},\check{v} \right\rangle_{\hat{V}',\hat{V}} =(f,v)_{D} + \alpha(|\Gamma|\hat{\psi},\check{v} )_{\Lambda}  \ \ \forall v \in V_D,\label{var1}\\
&&\Big(\tilde{K}|\Sigma|\cfrac{d\hat{u}}{ds},\cfrac{d\hat{v}}{ds}\Big)_{\Lambda}+\hat{\alpha}(|\Gamma|\hat{u},\hat{v})_{\Lambda}+\left\langle |\Gamma|\overline{\phi},\hat{v} \right\rangle_{\hat{V}',\hat{V}}=(|\Sigma|\overline{\overline{g}},\hat{v})_{\Lambda}+\hat{\alpha}(|\Gamma|\hat{\psi},\hat{v})_{\Lambda}\ \ \forall \hat{v} \in \hat{V} \label{var2},
\end{eqnarray}
being $\overline{\overline{g}}(s)=\frac{1}{|\Sigma(s)|}\int_{\Sigma(s)}g~d\sigma$. 
The terms multiplied by coefficients $\alpha$ and $\hat{\alpha}$ in the constraint equations \eqref{var1}-\eqref{var2} represent a consistent correction, as at the minimum $\check{u}=\hat{u}=\hat{\psi}$. However this correction allows to have well posed problems on each sub-domain independently from the prescribed boundary conditions, provided that $\alpha,\hat{\alpha}>0$. This is particularly relevant, as one of the key advantages of the proposed approach is to provide a methodology ready for domain decomposition on non conforming meshes. 
In addition, the discrete problem deriving from the optimization formulation is well posed without requiring the introduction of complex stabilization terms. The above formulation can be extended to accommodate multiple intersecting segments and different couplings between the 3D and 1D domain. The interested reader is referred to the previous works on the subject for further details \cite{BGS3D1D2022,BGS3D1Ddisc,BGSV3D1D2022}.

\section{Discrete problem} \label{discr}
Let us briefly recall here the discrete formulation of problem~\eqref{minJ}-\eqref{var2}, in the simplified case of a single inclusion. 
It is to remark that this is formally identical to the one already described in \cite{BGS3D1D2022}, also for the general case of multiple intersecting 1D domains. 
Indeed, the focus of the present work is on the application of the XFEM, which does not affect the structure of the discrete system. The choice of the enrichment function and of the quadrature formulas, which are instead the main novelty content of this work, are thoroughly discussed in the next sections. 

As the inclusion is reduced to the centerline $\Lambda$, we extend the domain $D$ to cover to the whole $\Omega$. Then we build a mesh $\mathcal{T}$ on $\Omega$ made of $N_\tau$ tetrahedral elements $\tau_j$, i.e. $\mathcal{T}=\bigcup_{j=1}^{N_\tau} \tau_j$, whose position in space is independent from the position of the 1D domain $\Lambda$. On this mesh we choose a set of finite element basis functions $\left\lbrace \varphi_i\right\rbrace_{i=1,\ldots,N}$, such that the discrete counterpart of unknown $u$ is $U=\sum_{i=1}^N U_i \varphi_i$. 
We proceed similarly for variables $\hat{u}$, $\bar{\phi}$ and $\hat{\psi}$, by first defining on $\Lambda$ three independent meshes and basis functions sets: mesh $\hat{\mathcal{T}}$ and functions $\{\hat{\varphi}_i\}_{i=1,\ldots,\hat{N}}$ for $\hat{u}$, mesh $\tau^\phi$ and functions $\{\theta_i\}_{i=1,\ldots,N^\phi}$ for $\bar{\phi}$, and mesh $\tau^\psi$ and functions $\{\eta_i\}_{i=1,\ldots,N^\psi}$ for $\hat{\psi}$, ending up with the following discrete counterparts for the three variables, defined respectively as:
\begin{displaymath}
\hat{U}=\sum_{i=1}^{\hat{N}}\hat{U}_{i}\hat{\varphi}_{i}, \quad \Phi=\sum_{i=1}^{N^{\phi}}\Phi_{i}\theta_{i}, \quad \Psi=\sum_{i=1}^{N^{\psi}}\Psi_{i}\eta_{i}.
\end{displaymath}
The discrete problem is obtained by replacing the above definitions in equations~\eqref{var1}-\eqref{var2}. The discrete functional is then defined as follows:
\begin{displaymath}
J_\delta=\frac12 \left(\|U_{|_\Lambda}-\Psi\|_{L^2(\Lambda)}^2+\|\hat{U}-\Psi\|_{L^2(\Lambda)}^2\right),
\end{displaymath}
i.e. exploiting the regularity of the discrete variables to directly compute the restriction on $\Lambda$ of $U$ and using the $L^2$-norm to compute the coupling mismatch. 
Then we collect the integrals of the basis functions into the matrices:
\begin{align*}
&\bm{A} \in \mathbb{R}^{N\times N} \text{ s.t. } (A)_{kl}=\int_{\Omega}\bm{{K}}\nabla\varphi_k\nabla\varphi_l ~d\omega+\alpha \int_{\Lambda}|\Gamma(s)|{\varphi_k}_{|_{\Lambda}}{\varphi_l}_{|_{\Lambda}} ds,\\[0.8em]
&\bm{\hat{A}} \in \mathbb{R}^{\hat{N}\times \hat{N}} \text{ s.t. } (\hat{A})_{kl}=\int_{\Lambda}\bm{\tilde{K}}|\Sigma(s)|\frac{d\hat{\varphi}_{k}}{ds}\frac{d\hat{\varphi}_{l}}{ds} ~ds+\hat{\alpha}\int_{\Lambda}|\Gamma(s)|\hat{\varphi}_{k}\hat{\varphi}_{l}~ds,
\end{align*}
\begin{align*}
&\bm{B} \in \mathbb{R}^{N\times N^{\phi}} \text{ s.t. } (B)_{kl}=\int_{\Lambda}|\Gamma(s)|{{\varphi_{k}}_{|_{\Lambda}}\theta_{l}}~ds,\\
&\bm{\hat{B}} \in \mathbb{R}^{\hat{N}\times N^{\phi}} \text{ s.t. } (\hat{B})_{kl}=\int_{\Lambda}|\Gamma(s)|{\hat{\varphi}_{k}~\theta_{l}}~ds,\\
&\bm{C}^{\alpha} \in \mathbb{R}^{N\times N^{\psi}} \text{ s.t. } (C^{\alpha})_{kl}=\alpha\int_{\Lambda}|\Gamma(s)|{\varphi_k}_{|_{\Lambda}}\eta_{l}~ds,\\
&\bm{\hat{C}}^{\alpha} \in \mathbb{R}^{\hat{N}\times N^{\psi}} \text{ s.t. } (\hat{C}^{\alpha})_{kl}=\hat{\alpha}\int_{\Lambda}|\Gamma(s)|{\hat{\varphi}}_{k}~\eta_{i,l}~ds,
\end{align*}
and vectors 
\begin{align*}
&f\in \mathbb{R}^N \text{ s.t. } f_k=\int_{\Omega}f\varphi_k~d\omega,\\
&g\in \mathbb{R}^{\hat{N}} \text{ s.t. } (g)_k=\int_{\Lambda}|\Sigma(s)|\overline{\overline{g}}~\hat{\varphi}_{k}~ds.
\end{align*}
We end up in the following form of the constraints:
\begin{align}
&\bm{A}U-\bm{B}\Phi-\bm{C}^{\alpha}\Psi=f\label{eq1discr},\\
&\bm{\hat{A}}\hat{U}+\bm{\hat{B}}\Phi-\bm{\hat{C}}^{\alpha}\Psi=g,\label{eq2discr}
\end{align}
where, with a notation overload, we denoted the array of degrees of freedom with the same symbol of the corresponding discrete function. 
We proceed similarly for the functional, which, after defining:
\begin{align*}
&\bm{G} \in \mathbb{R}^{N \times N} \text{ s.t. } (G)_{kl}=\int_{\Lambda}{\varphi_k}_{|_{\Lambda}}{\varphi_l}_{|_{\Lambda}}ds,\\
&\bm{\hat{G}} \in \mathbb{R}^{\hat{N} \times \hat{N}} \text{ s.t. } (\hat{G})_{kl}=\int_{\Lambda}{\hat{\varphi}}_{k}~{\hat{\varphi}}_{l}~ds,\\
&\bm{G^{\psi}} \in \mathbb{R}^{N^{\psi} \times N^{\psi}}  \text{ s.t. } (G^{\psi})_{kl}=\int_{\Lambda}\eta_{k}~\eta_{l}~ds,\\
&\bm{C} \in \mathbb{R}^{N\times N^{\psi}} \text{ s.t. } (C)_{kl}=\int_{\Lambda}{\varphi_k}_{|_{\Lambda}}\eta_{l}~ds,\\
&\bm{\hat{C}} \in \mathbb{R}^{\hat{N}_i\times N^{\psi}} \text{ s.t. } (\hat{C})_{kl}=\int_{\Lambda}{\hat{\varphi}}_{k}~\eta_{l}~ds,
\end{align*}
reads:
\begin{displaymath}
J_\delta=\frac12\left( U^T\bm{G}U-U^T\bm{C}\Psi-\Psi^T\bm{C}^TU+\hat{U}^T\bm{\hat{G}}\hat{U}-\hat{U}^T\bm{\hat{C}}\Psi-\Psi^T\bm{\hat{C}}^T\hat{U}+2\Psi^T\bm{G^{\psi}}\Psi\right).
\end{displaymath}
The discrete problem:
\begin{eqnarray*}
	&&\min_{(\Phi,\Psi)}J_\delta(\Phi,\Psi),\\
	&& \text{ subject to  \eqref{eq1discr}-\eqref{eq2discr}} 
\end{eqnarray*}
can be solved either by resorting to the corresponding saddle-point system of optimality condition and by means of a gradient based descent method applied to the unconstrained minimization problem obtained replacing the linear constraints into the functional. 
Further details are available in \cite{BGSV3D1D2022,BGS3D1Ddisc}.

\section{Application of the XFEM}\label{xfem_func}
In the definition of the discrete function $U$ we have denoted by $\{\varphi_k\}_{k=1,\ldots,N}$ a generic set of basis functions. We can now better define such basis functions, splitting them into two sets: the set of standard linear Lagrangian basis functions, denoted by $\{\varphi_k^s\}_{k\in \mathcal{I}}$, with $\mathcal{I}$ the set of their DOF indexes, and the set of the enrichment basis functions $\{\varphi_k^e\}_{k\in \mathcal{I^\star}}$, with $\mathcal{I^\star}$ the corresponding set of DOF indexes. The enrichment basis functions are built starting from a global enrichment function and then applying the partition of unity method.  

The global enrichment function needs to account for the the irregular behavior that is expected for the discrete solution. Our choice of enrichment function follows from the results in \cite{Gjerde2019}. 
We will denote by $\zeta(\bm{x})$ the global enrichment function, which is defined in a different way if the inclusion entirely crosses the domain of interest, or if, instead, it is embedded in the domain. In this latter case, indeed, the enrichment function also needs to control the shape of the solution around the endpoints of the inclusion. Let us start by considering a single inclusion and let us denote by $d_\ell(\bm{x})$ the distance of a generic point $\bm{x}$ from the line $\ell$ passing through $\Lambda$ (see again Figure~\ref{fig:notation}). In case $\Lambda$ crosses the domain from side to side, we define $\zeta$ as:
\begin{equation}
\zeta(\bm{x})=\zeta^\natural(\bm{x}):=\begin{cases}
-\log(d_\ell(\bm{x})) &\text{if } d_\ell(\bm{x})>R\\
-\log(R) &\text{if } d_\ell(\bm{x})\leq R,
\end{cases}\label{eq:zeta_bequadro}
\end{equation}
\begin{figure}
	\centering
	\begin{subfigure}[t]{.45\textwidth}
		\centering
		\includegraphics[width=1.05\linewidth]{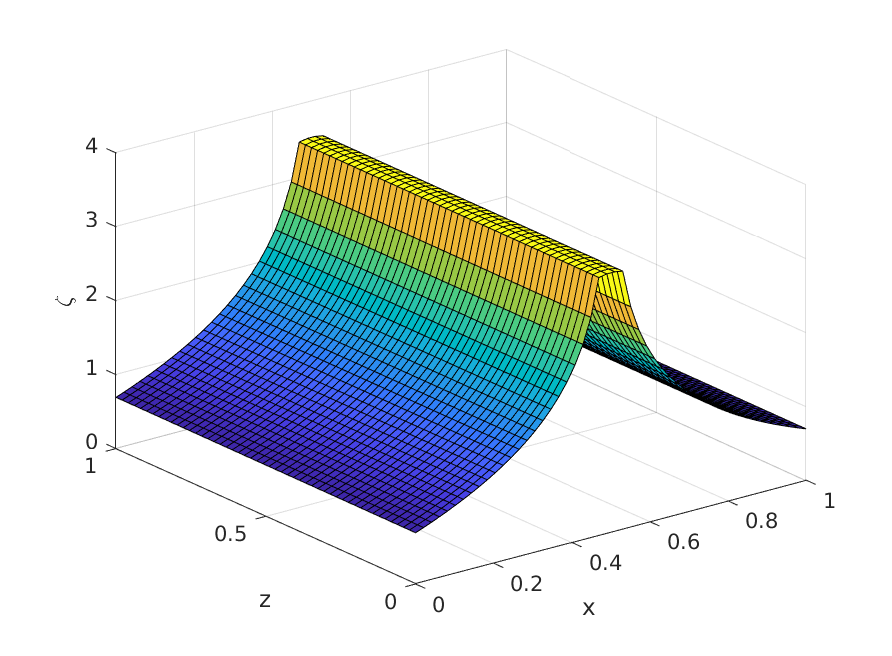}
		\label{}
		\caption{$\zeta=\zeta^\natural$}
	\end{subfigure}\hfill
	\begin{subfigure}[t]{.45\textwidth}
		\centering
		\includegraphics[width=1.05\linewidth]{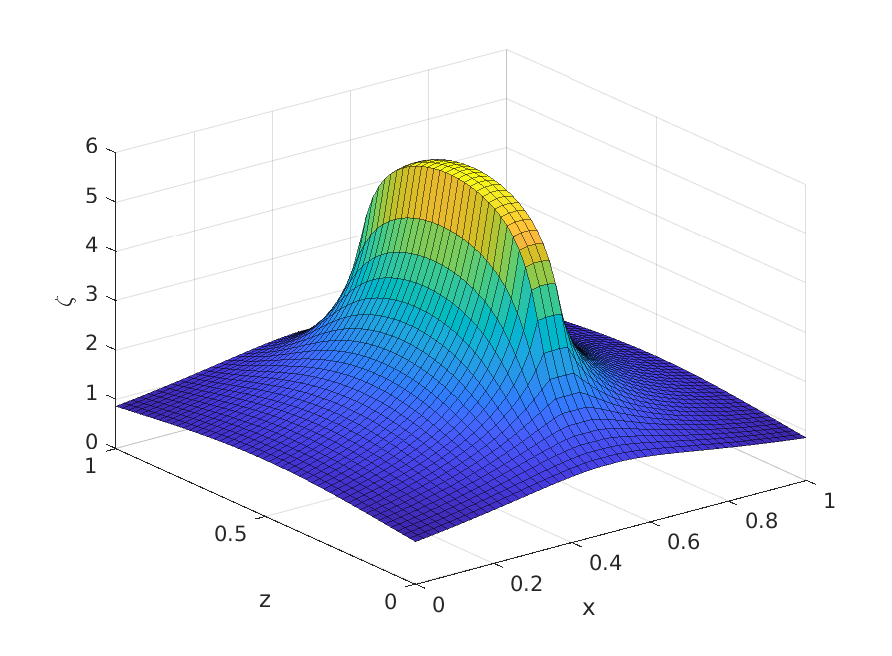}
		\label{}
		\caption{$\zeta=\zeta^\flat$}
	\end{subfigure}
	\caption{Enrichment function $\zeta$ on a plane containing inclusion centreline. $\Omega=[0,1]^3$, $\Lambda$ aligned with the $z$-axis, $R=10^{-2}$. For $\zeta=\zeta^\flat$, $\Lambda$ extending from $z=0.2$ to $z=0.8$.}
	\label{fig:zeta_general}
\end{figure}
being $R$, as before, the radius of the original 3D inclusion $\Sigma$.

If, instead, the endpoints $\bm{x}_0$ and $\bm{x}_S$ of $\Lambda$ lie inside $\Omega$, we define
\begin{equation}
\zeta(\bm{x})=\zeta^\flat(\bm{x})=\begin{cases}
\log{\left(\cfrac{||\bm{x}-\bm{x}_S||+L+\bm{\tau_\Lambda}\cdot(\bm{x}_0-\bm{x})}{||\bm{x}-\bm{x}_0||+\bm{\tau_\Lambda}\cdot (\bm{x}_0-\bm{x})}\right)} &\text{if } d_\ell(\bm{x})>R,\smallskip\\
\log{\left(\cfrac{||\bm{x}_R-\bm{x}_S||+L+\bm{\tau_\Lambda}\cdot(\bm{x}_0-\bm{x}_R)}{||\bm{x}_R-\bm{x}_0||+\bm{\tau_\Lambda}\cdot (\bm{x}_0-\bm{x}_R)}\right)} &\text{if }  d_\ell(\bm{x})\leq R
\end{cases}\label{eq:zeta_bemolle}
\end{equation}
where $\bm{\tau_\Lambda}$ is the unit tangent vector to $\Lambda$, $\bm{x}_R$ is the projection of $\bm{x}$ on the infinite cylindrical surface $\mathcal{G}$ containing $\Gamma$, $\bm{x}_0, \ \bm{x}_S$ are the endpoints of $\Lambda$ and $L=||\bm{x}_S-\bm{x}_0||$.
Let us observe that
\begin{equation}
\lim_{L\rightarrow \infty}\frac{1}{4\pi}\zeta^\flat(\bm{x})\approx -\frac{1}{2\pi}\zeta^\natural(\bm{x})
\end{equation}
is a relation usually used in electromagnetism to approximate the potential of an infinite length line charge.
Functions $\zeta^\natural(\bm{x})$ and $\zeta^\flat(\bm{x})$ are shown in Figure~\ref{fig:zeta_general} on a plane containing $\Lambda$.

Let us now consider a cylinder $\Delta$ in $\Omega$ with centreline coinciding with $\Lambda$ and having constant cross-section radius $\rho\geq R$ and 
let us further denote by $\mathcal{T}_\Delta$ the subset of mesh elements in $\mathcal{T}$ having an intersection with the cylinder $\Delta$ of non null measure, i.e. $\mathcal{T}_\Delta:=\left\lbrace \tau \in \mathcal{T}: |\tau\cap\Delta|>0\right\rbrace$. We denote by $\mathcal{J}\subset\mathcal{I}$ the degree of freedom indexes $k\in \mathcal{I}$ such that the support of the standard basis function $\varphi^s_k$ has a non empty overlap with an element in $\mathcal{T}_\Delta$, i.e. $\mathcal{J}:=\left\lbrace k\in\mathcal{I}: \exists \tau_j\in \mathcal{T}_\Delta, \text{supp}(\varphi^s_k)\cap \tau_j \neq \emptyset \right\rbrace$. 
We also introduce a continuous ramp function $r_{\mathcal{J}}(\bm{x})$, equal to one inside $\mathcal{T}_\Delta$ and linearly vanishing to zero outside $\mathcal{T}_\Delta$, obtained as $r_{\mathcal{J}}(\bm{x}):=\sum_{k\in\mathcal{J}} \varphi^s_k(\bm{x})$.
For $k\in \mathcal{J}$ we then define $\bar{\varphi}_k^e(\bm{x}) =\varphi^s_k(\bm{x}) \zeta(\bm{x})  r_{\mathcal{J}}(\bm{x})$ and finally $\varphi_k^e(\bm{x})=\bar{\varphi}_k^e(\bm{x})-\bar{\varphi}_k^e(\bm{x}_k)$, such that the enrichment basis functions are zero-valued in the mesh vertexes $\bm{x}_k$.
Following the XFEM paradigm \cite{XFEMreview}, the effect of the enrichment is local, in a neighborhood of the 1D domain $\Lambda$, depending on the the chosen value of $\rho$. 

The extension to the case of multiple inclusions is quite straightforward, by simply using the superposition effect. Let us consider $\mathcal{L}$ inclusions $\Lambda_i$, and, for each inclusion, let us define a cylinder $\Delta_i$ with a centreline coinciding with $\Lambda_i$ and radius $\rho_i > R_i$, being $R_i$ the radius of the 3D inclusion $\Sigma_i$. 
We then define
$$\mathcal{T}_{\Delta}^i:=\left\lbrace \tau \in \mathcal{T}: |\tau\cap\Delta_i|>0\right\rbrace$$ and 
$$\mathcal{J}_i:=\left\lbrace k\in\mathcal{I}: \exists \tau_j\in \mathcal{T}_{\Delta}^i, \text{supp}(\varphi^s_k)\cap \tau_j \neq \emptyset \right\rbrace.$$
A different enrichment function is defined for each inclusion, namely
\begin{displaymath}
\varphi^{e, i}_k(\bm{x}) =\varphi^s_k(\bm{x}) \zeta_i(\bm{x})  r_{\mathcal{J}}^i(\bm{x}),
\end{displaymath}
with $r_{\mathcal{J}}^i(\bm{x}):=\sum_{k\in\mathcal{J}_i} \varphi^s_k(\bm{x})$ and $\zeta_i(\bm{x})$ defined as in \eqref{eq:zeta_bequadro} or \eqref{eq:zeta_bemolle} depending on $\Lambda_i$. The discrete approximation of the unknown $u$ is then defined as:
$$U=\sum_{k\in \mathcal{I}} U_k^s\varphi^s_k(\bm{x}) + \sum_{i=1}^{\mathcal{L}}\sum_{k\in \mathcal{J}_i} U_{k}^{e,i}\left(\varphi^{e,i}_k(\bm{x})\right).$$
In practice, the unknowns and the corresponding basis functions are numbered consecutively, giving $N$ total unknowns.

Please note that the case of intersecting inclusions is contained in the above presentation, as we can simply split the intersecting centrelines into sub-segments meeting in one of their endpoints.  

\section{Numerical integration} \label{quad}
\begin{figure}
	\centering
	\includegraphics[width=0.75\textwidth]{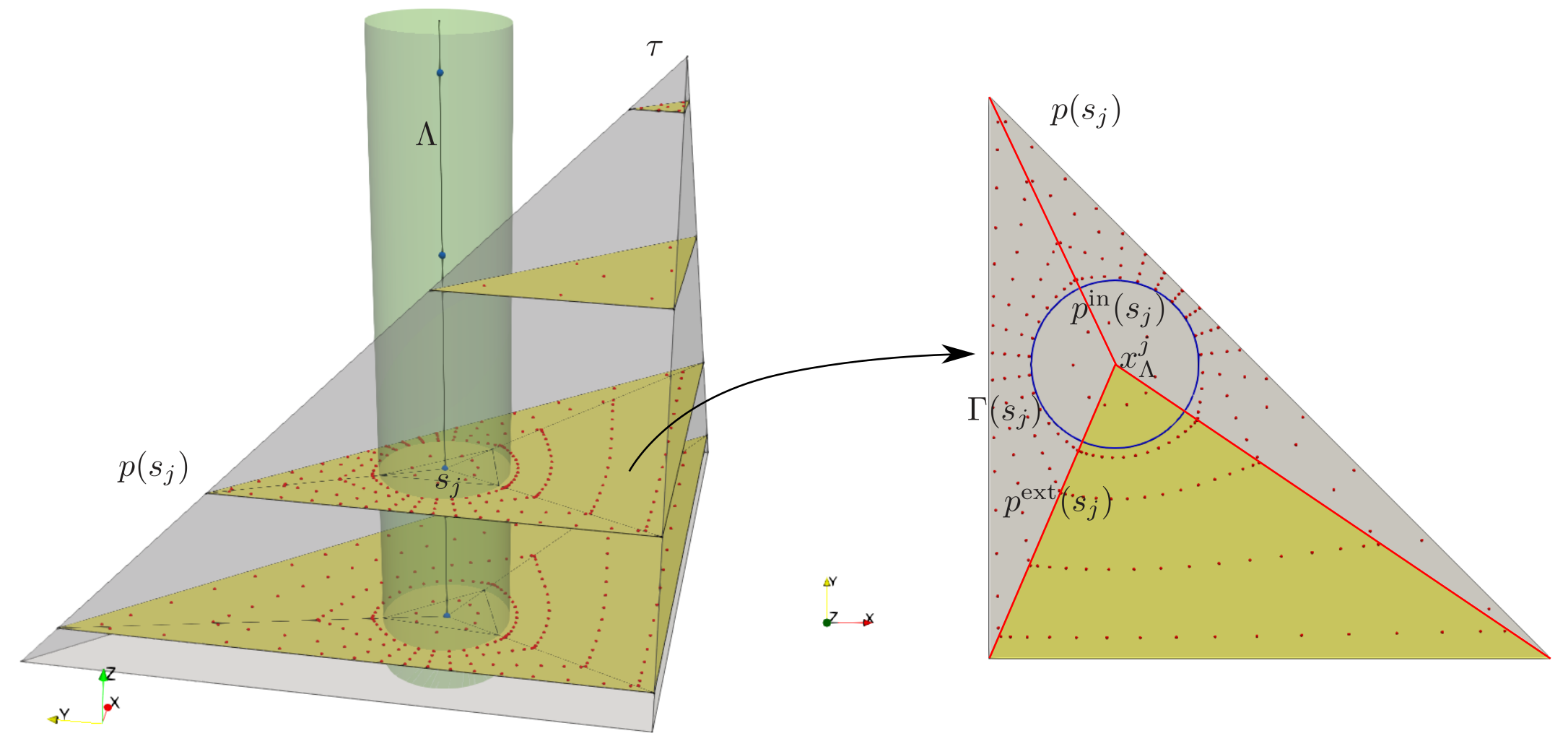}
	\caption{Description of numerical integration strategy}
	\label{Fig_quad}
\end{figure}
\begin{table}
	\centering
	\caption{Numerical quadrature errors for different numbers of integration points}
	\label{Tab_quad}
	\begin{tabular}{ccc|c|c|c}
		$n_\Lambda$ & $n_\mathfrak{r}$ & $n_\theta$ & $N_{\text{pt}}$ & error - $R=0.1$ & error - $R=0.3$\\
		\hline
		1& 3&  5&  33& 6.94e-05& 4.28e-06\\
		1& 4&  7&  59& 9.85e-08& 2.97e-09\\
		1& 6&  9& 111& 6.45e-12& 1.75e-12\\ 
		1& 8& 12& 195& 4.57e-16& 1.67e-16\\
	\end{tabular}
\end{table}
A key aspect for the successful application of the XFEM lies in the numerical quadrature of the enrichment basis functions. 
Given the irregular behavior of such functions, customized strategies need to be adapted, often relying on a sub-division of the three-dimensional domain conforming to the interfaces.
Here, the devised approach exploits the known behavior of function $\zeta(\bm{x})$ and is capable of correctly capturing the curvilinear boundary of the interface.
Let us start with the case of an isolated inclusion. With reference Figure~\ref{Fig_quad}, let us consider a tetrahedron $\tau\in\mathcal{T}$, intersected by one inclusion $\Lambda$ with radius $R$. 
Let us denote by $s_{v_0}\geq s_{v_1}\geq s_{v_2}\geq s_{v_3}$ the curvilinear abscissas of the projections on $\Lambda$ of the four vertexes of $\tau$. We remark that it is possible that some of these projection points coincide, when $\Lambda$ is orthogonal to one of the faces of the element, as it is the case of Figure~\ref{Fig_quad}.
Then, considering a generic enrichment $\varphi^e$ we have:
\begin{displaymath}
\int_{\tau} \varphi^e(\bm{x}) \text{d}\bm{x} = \sum_{t = 0}^2 \int_{s_{v_t}}^{s_{v_{t + 1}}} \left( \int_{p(s)} \varphi^e(\bm{x}) \text{d}\sigma \right) := \sum_{t = 0}^2 \int_{s_{v_t}}^{s_{v_{t + 1}}} f_\zeta(s) \text{d}s,
\end{displaymath}
in which $p(s)$ is the polygonal region given by the intersection of $\tau$ with a plane orthogonal to $\Lambda$ at $s\in [s_{v_0}, s_{v_3}]$. In each interval $\Lambda^t:=[s_{v_t}, s_{v_{t+1}}]$, $t=0,\ldots,2$, the function representing the surface area of $p(s)$ is smooth, and consequently $f_\zeta(s)$ is smooth.
A Gaussian 1D quadrature rule with $n_\Lambda$ nodes can be efficiently adopted to integrate $f_\zeta(s)$ in each $\Lambda^t$, requiring the computation of values $f_\zeta(s_j)$ at integration nodes $s_j\in (s_{v_t}, s_{v_{t+1}})$, $j=1,\ldots,n_\Lambda$.
The strategy to compute integrals $f_\zeta(s_j)=\int_{p(s_j)} \zeta(\bm{x})$ on the regions $p(s_j)$, instead, is different depending on the position of $p(s_j)$. Indeed, if $p(s_j)$ does not contain the irregularity interface of function $\zeta$, we adopt standard quadrature. This is the case, for example, of the two top triangular regions in Figure~\ref{Fig_quad}. Whereas, when $p(s_j)$ contains the interface, as in the two bottom triangular regions in Figure~\ref{Fig_quad}, the integration is performed combining the approaches proposed in \cite{FALLETTA2015106} and \cite{MONEGATO1999201}, adapted to the present case, and described in the following. 
Let us denote by $\bm{x}^j_\Lambda$ the point at the intersection between $\Lambda$ and the plane containing $p(s_j)$, and let us denote by $\Gamma(s_j)$ the irregularity interface of $\zeta$ at $s_j$. We remark that this actually coincides with the intersection of $p(s_j)$ with the lateral surface of the original 3D inclusion, see Figure~\ref{Fig_quad}, right.
Furthermore, let us call $p^\text{in}(s_j)$ the portion of $p(s_j)$ inside $\Gamma(s_j)$, and $p^\text{ext}(s_j)$ the portion outside $\Gamma(s_j)$. The regions $p^\text{in}(s_j)$ and $p^\text{ext}(s_j)$ are each covered by triangular regions with one vertex in $\bm{x}^j_\Lambda$, as illustrated in Figure~\ref{Fig_quad}, right. In this case we have that $\Gamma(s_j)$ is entirely contained in $p(s_j)$, but such covering can also be determined when it is only partially contained. More details are available in \cite{FALLETTA2015106}.
Now we apply a first mapping $t:[x,y]\mapsto [x^\star,y^\star]$ from each triangular region to the reference triangle, with $\bm{x}^j_\Lambda$ being mapped to the origin of the reference frame $(x^\star,y^\star)$. We remark that, in such reference frame, the portion of $\Gamma(s_j)$ contained in the triangular region is mapped to an ellipse, centered in the origin. Then we apply a rotation $\varrho:[x^\star,y^\star]\mapsto [\tilde{x},\tilde{y}]$ to align the axis $x^\star y^\star$  with the principal axis of this ellipse, whose equation is $\frac{\tilde{x}^2}{\lambda_1^2}+\frac{\tilde{y}^2}{\lambda_2^2}=R^2$ in the new frame.  Finally we introduce a polar transformation $\Upsilon:[\mathfrak{r},\theta]\mapsto [\tilde{x},\tilde{y}]$, depending on a parameter $q$ and defined as:
\begin{equation*}
\begin{cases}
\tilde{x} = \lambda_1 R \mathfrak{r}^q \cos{\theta}\\
\tilde{y} = \lambda_2 R \mathfrak{r}^q \sin{\theta}.
\end{cases}
\end{equation*} 
Now we choose $n_\mathfrak{r}$ Gaussian quadrature nodes along $\mathfrak{r}$ and $n_\theta$ nodes along $\theta$ which are then mapped back to the physical reference frame $(x,y)$. The three changes of variables allow to correctly integrate the enrichment function close to the curvilinear interface $\Gamma$, since the value of $q$ can be chosen to obtain a clustering of the nodes towards the border of the ellipse, where the function has a steep gradient. Higher values of $q$ correspond to a higher clustering. A value $q=1$ is used for the regions inside $\Gamma$, where $\zeta$ is constant, whereas a value $q=3$ is employed for the regions outside $\Sigma$. We remark that the above quadrature strategy also applies to general polyhedrons.

As an example, we integrate the  function $\zeta^\natural(d_\Lambda(\bm{x}))$ defined in\eqref{eq:zeta_bequadro} over a unit edge cubic domain~$[0, 1]^3$, where $\Lambda$ coincides with the vertical edge of the cube passing through the origin.
We chose this simple geometry to allow for the computation of the exact integral.
Two values of $R$ are proposed: $R=0.1$ and $R=0.3$. 
The obtained results are reported Table~\ref{Tab_quad}, showing that with the proposed strategy, it is possible to compute the integral from single precision up to machine precision. 
In the table, $N_{\text{pt}}$ represents the total number of quadrature points.
The values in columns $n_\mathfrak{r}$ and $n_\theta$ refer to the number of quadrature nodes selected in the external regions $p^\text{ext}(\cdot)$.
Constant values of $n_\mathfrak{r} = 1$ and $n_\theta = 1$ are used for the internal regions $p^\text{in}(\cdot)$.
Moreover, in this particular case, a single node along $\Lambda$ is sufficient, considering the simple geometry of the domain and the regularity of the integrated function $f_{\zeta^\natural}(s)$. 


%
\begin{figure}
	\centering
	\begin{subfigure}[t]{.225\textwidth}
		\centering
		\includegraphics[width=1\linewidth]{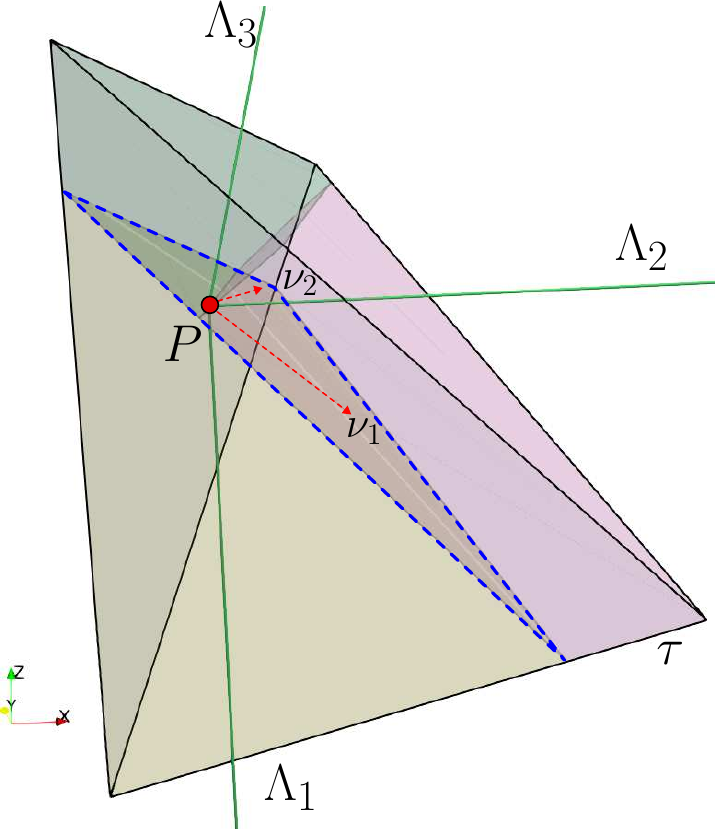}
		\caption{Tetrahedron, split selecting $\Lambda_1$ and $\Lambda_2$.}
		\label{fig:splitTetra-a}
	\end{subfigure}
	\begin{subfigure}[t]{.26\textwidth}
		\centering
		\includegraphics[width=1\linewidth]{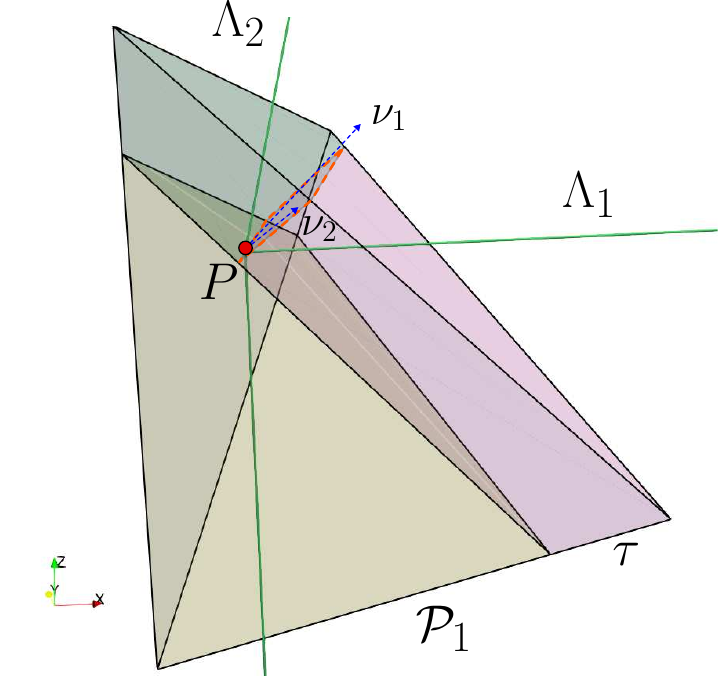}
		\caption{Further split selecting $\Lambda_1$ and $\Lambda_2$.}
		\label{fig:splitTetra-b}
	\end{subfigure}
	\begin{subfigure}[t]{.21\textwidth}
		\centering
		\includegraphics[width=1\linewidth]{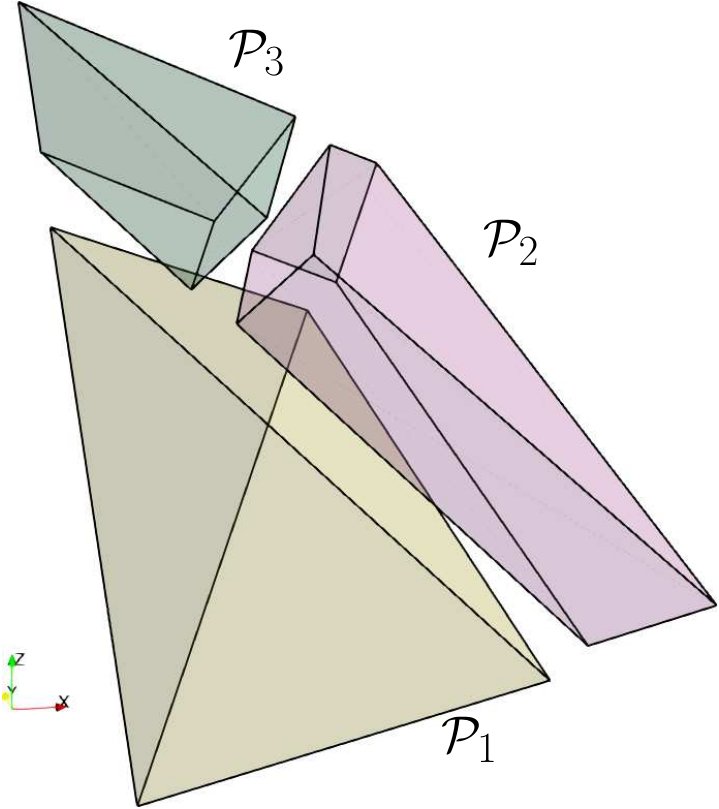}
		\caption{Tetrahedron split.}
		\label{fig:splitTetra-c}
	\end{subfigure}
	\begin{subfigure}[t]{.22\textwidth}
		\centering
		\includegraphics[width=1\linewidth]{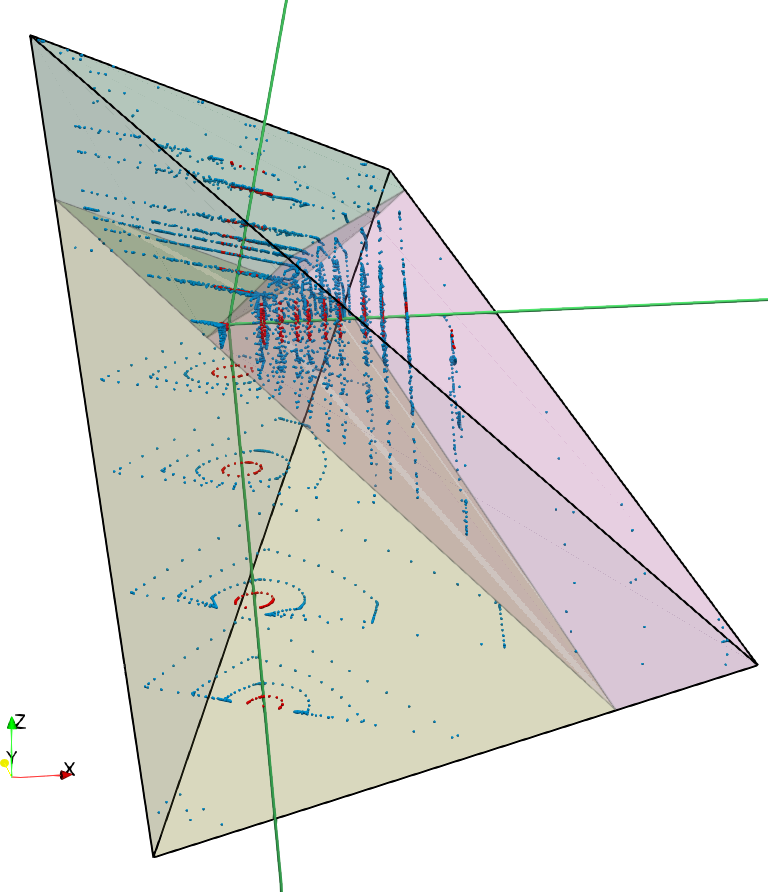}
		\caption{Quadrature points generated; red points are inside the inclusions.}
		\label{fig:splitTetra-d}
	\end{subfigure}
	\caption{Description of the split strategy.}
	\label{fig:splitTetra}
\end{figure}
%

The above procedure is generalized to the case of multiple inclusions as follows. 
If multiple non intersecting inclusions pass through a single tetrahedron $\tau \in \mathcal{T}$, or if a single segment ends within a tetrahedron, it is sufficient to split the element into sub-cells such that each sub-cell only contains up to one inclusion, entirely crossing it.
Then, we use the quadrature strategy proposed above in cells containing an inclusion, or a classic one, if the considered sub-cell contains no segments. 

The case of multiple segments intersecting in point $P$ in a tetrahedral cell $\tau \in \mathcal{T}$ also requires a splitting into sub-cells containing up to a single inclusion (or a portion of a single inclusion). 
We remark that, in this case, the enrichment function $\zeta=\zeta_i^{\flat}$ needs to be used, as segment endpoints (at least those matching with $P$) are inside $\Omega$. 

We choose to split elements with the following strategy. Let us consider $\mathcal{L}\geq 2$ intersecting segments in $\tau \in \mathcal{T}$, locally numbered as $\Lambda_i$, $i=1,\ldots,\mathcal{L}$. We select directions $\bm{\nu}_1$ and $\bm{\nu}_2$ as the sum and the external product of the unit tangent vectors of the first two centrelines, respectively, i.e. $\bm{\nu}_1=\bm{\tau}_{\bm{\Lambda}_1}+\bm{\tau}_{\bm{\Lambda}_2}$, and $\bm{\nu}_2=\bm{\tau}_{\bm{\Lambda}_1}\wedge\bm{\tau}_{\bm{\Lambda}_2}$. Then we cut cell $\tau$ along the plane containing $\bm{\nu}_1$ and $\bm{\nu}_2$ and passing through $P$. This generates two sub-cells, and the procedure is replicated on each sub-cell. If the sub-cell contains more than two inclusions it is split again along a cutting direction, chosen as above and depending on the local (arbitrary) renumbering of the inclusions in the sub-cell itself. If instead a sub-cell contains one or no inclusions, it is left unchanged. The process is recursively applied to each newly generated sub-cell until all sub-cells contain less than $2$ inclusions.
An example is shown in Figure~\ref{fig:splitTetra}, for a cell containing three inclusion. The first cut is performed along the plane containing  $\bm{\nu}_1$ and $\bm{\nu}_2$ in Figure~\ref{fig:splitTetra-a} and passing through $P$, thus generating two sub-cells $\mathcal{P}_1$ and $\mathcal{P}_2$. Cell $\mathcal{P}_1$ contains a single inclusion, and requires no further splitting. Cell $\mathcal{P}_2$ instead still contains $2$ inclusions, locally renumbered as $\Lambda_1$ and $\Lambda_2$, see Figure~\ref{fig:splitTetra-b}, and thus it is further split along the plane containing $\bm{\nu}_1$ and $\bm{\nu}_2$ passing through $P$, giving cells $\mathcal{P}_2$ and $\mathcal{P}_3$, see Figure~\ref{fig:splitTetra-c}.
Finally, on each sub-cell we apply the quadrature strategy outlined in Section~\ref{quad}. Figure~\ref{fig:splitTetra-d} shows quadrature nodes for the proposed example.

\section{Numerical results} \label{num_res}
The following section is devoted to the presentation of five numerical tests in order to validate and show the effectiveness of the proposed approach. In the following we will denote by $N$ the number of degrees of freedom for variable $U$, which for a fixed  mesh can vary according to the radius $\rho$ of the enrichment cylinder $\Delta$. Let us remark that the case $\rho>0$ corresponds to the optimization based domain decomposition method with the use of the extended finite elements for the 3D variable, while for $\rho=0$ we end up in the same optimization based approach but with standard finite elements for the 3D variable. For what concerns the 1D variables, piecewise linear continuous basis functions are used for $\hat{U}$ on an equally spaced mesh $\hat{\mathcal{T}}$ and for $\Psi$ on an equally spaced mesh $\mathcal{T}^\psi$, whereas piecewise constant basis functions are used for $\Phi$ on an equally spaced mesh $\mathcal{T}^\phi$. The refinement level of the 1D meshes $\hat{\mathcal{T}}$, $\mathcal{T}^\psi$ and $\mathcal{T}^\phi$ is related to the refinement level of the 3D mesh $\mathcal{T}$ and not to the number of degrees of freedom $N$. More in details, denoted by $N_I$ the number of intersection points between an inclusion $\Lambda$ and the boundary of the elements in $\mathcal{T}$, mesh $\hat{\mathcal{T}}$ will have $2N_I$ nodes whereas meshes $\mathcal{T}^\psi$ and $\mathcal{T}^\phi$ will count $N_I/2$ nodes. The same is used for each segment in the case of multiple inclusions.  The analysis on the behavior of the method with respect to different refinement levels of the various meshes is available in the references, see \cite{BGS3D1D2022,BGS3D1Ddisc}.


\subsection{3D problem with singular source term} 
\label{test1}
\begin{figure}
	\centering
	\includegraphics[width=0.63\linewidth]{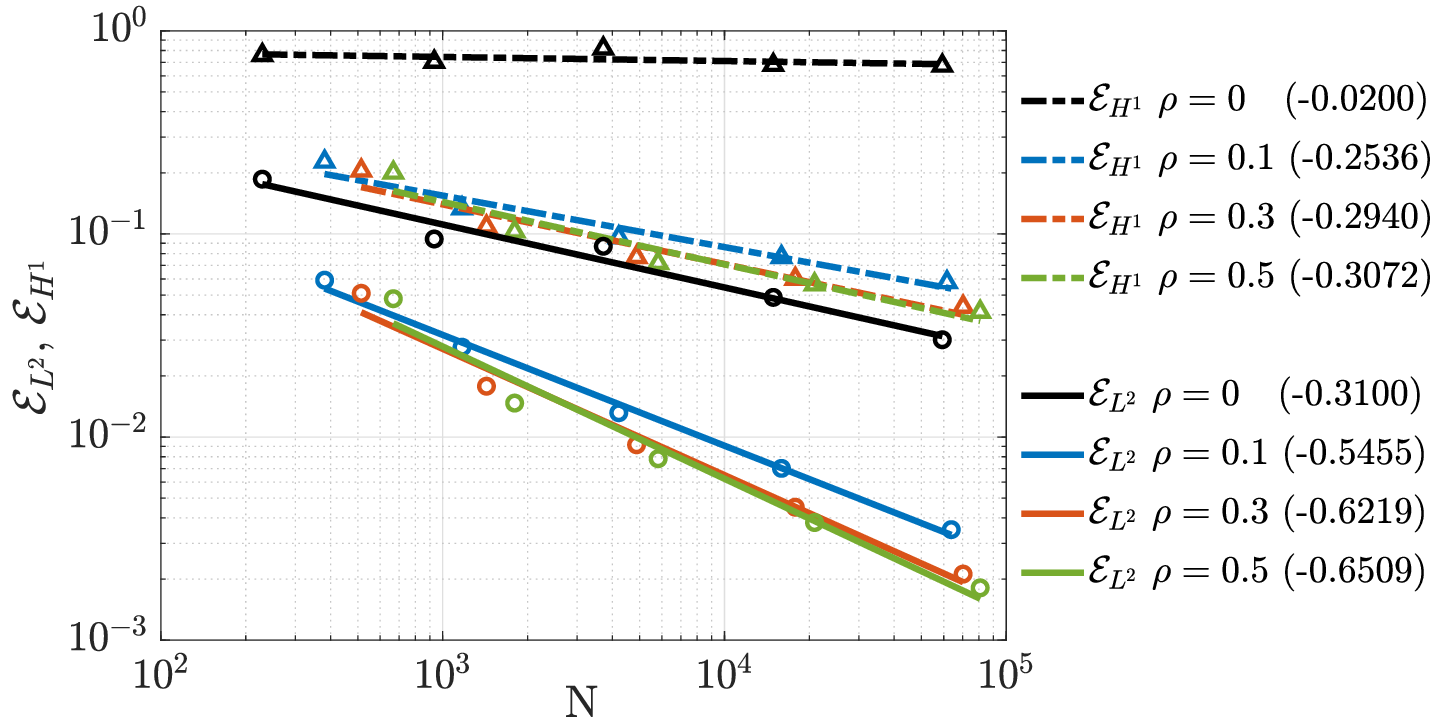}%
	\caption{Test~\ref{test1}, trend of the relative errors under mesh refinement. Dashed lines: relative $H^1$-norm of the error; full lines: relative $L^2$-norm of the error.}
	\label{fig:Test3D}
\end{figure}
The first numerical example concerns a 3D problem with a singular source term, and is used to validate the proposed XFEM setting through the comparison with a known analytical solution. Here, we will not solve a coupled 3D-1D problem, but a 3D problem with a known source term on a line. Consequently there is no need of using the optimization based coupling strategy. This example is therefore useful to investigate the effectiveness of the enrichment function shown in Section~\ref{xfem_func} in describing the expected behavior of the solution on coarse meshes and of the quadrature strategy described in Section~\ref{quad}.  

The test considers a cubic domain $\Omega=(-1,1)^3$ with a cylindrical inclusion $\Sigma=\lbrace(x,y,z): \sqrt{x^2+y^2}<R,~ z\in (-1,1)\rbrace,$ of radius $R=10^{-3}$. The inclusion is dimensionally reduced to a line and we numerically solve problem~\eqref{var1} with $\alpha=0$, $f=0$, $K=1$ and $\overline{\phi}=-\frac{1}{10\pi R}$. The obtained solution is compared to the analytical solution of the original equi-dimensional problem, chosen as:
\begin{displaymath}
u_{\text{ex}}=\begin{cases}
\frac{1}{10 \pi} \log(r) & \text{for} \ r>R \\
\frac{1}{10 \pi} \log(R) & \text{for} \ r\leq R
\end{cases}
\end{displaymath}
with $r=\sqrt{x^2+y^2}$.
We can observe that such analytical solution actually matches with the enrichment function $\zeta=\zeta^\natural$, but in the dimensionally reduced problem the flux $\overline{\phi}$ is placed at inclusion centreline and not at the inclusion boundary, thus introducing a modeling error. 

The problem is solved on five meshes with maximum element diameter ranging between $0.215$ and $0.034$. Homogeneous Neumann boundary conditions are prescribed on $\partial \Omega_{\mathrm{n}}=\lbrace (x,y,z): z=-1 \vee z=1\rbrace$, whereas Dirichelet boundary conditions, in accordance with the chosen exact solution, are set on $\partial \Omega_{\mathrm{d}}=\partial \Omega \setminus \overline{\partial \Omega_{\mathrm{n}}}$.
Convergence trends of the error between the computed and the analytical solution against the total number of degrees of freedom $N$ are reported in Figure~\ref{fig:Test3D} for the $L^2$ and $H^1$ relative norms. Four values of the enrichment area are considered, depending on the radius $\rho$ of cylinder $\Delta$ (see Section~\ref{xfem_func}): namely $\rho \in \{0,0.1,0.3,0.5\}$.
Table~\ref{tab:quad_node} reports the chosen quadrature parameters: $n_\Lambda$, $n_{\mathfrak{r}}$ and $n_\theta$ refer to the number of nodes along $\Lambda$, $\mathfrak{r}$ and $\theta$ respectively, as described in Section~\ref{quad}. For $n_{\mathfrak{r}}$ and $n_\theta$ we distinguish between the number of quadrature nodes used within the inclusion radius (\textit{in}), where the enrichment is constant, or outside (\textit{out}). Let us recall that the quadrature rule described in Section \ref{quad} is used only for elements which are intersected by $\Sigma$. The symbol $n_{\Delta}$, instead, denotes the number of nodes of a standard 3D Gaussian quadrature rule adopted on the tetrahedrons intersecting region $\Delta$ but not $\Sigma$: in these elements, indeed, we still need to integrate the enrichment functions, but, here, such functions have a continuous gradient.  
We also remark that the number of quadrature nodes used in a tetrahedron intersected by $\Sigma$ is typically larger than $n_\Lambda\times n_{\mathfrak{r}}\times n_\theta$, as it depends on the number of sub intervals $\Lambda^t$ (see Section~\ref{quad}) used and the number of sub-cells originated by the splitting. 
Further comments on this aspect are provided in the next examples, in which the same quadrature parameters reported in Table \ref{tab:quad_node} will be considered. 
Convergence trends in Figure~\ref{fig:Test3D} are close to the optimal ones for linear Lagrangian finite elements with regular data, and slightly improve if $\rho$ increases. The parameters reported in Table~\ref{tab:quad_node} are the proposed optimal choice: less nodes yield a decay in convergence trends and an upward shift of the error curves; more nodes lead to a small downward shift of the error curves.
If, instead, we choose $\rho=0$, i.e. we use standard finite elements, we have no convergence in the $H^1$ norm and a degraded $L^2$ convergence trend, which is in line with classical results for problems with singular data \cite{Scott1973}.

\begin{table}
	\renewcommand*{\arraystretch}{1.1}
	\centering
	\caption{Number of quadrature nodes}
	\label{tab:quad_node}
	\begin{tabular}{cccc}
		$n_\mathfrak{r}$ & $n_{\theta}$ & $n_{\Lambda}$ & $n_{\Delta}$ \\
		\hline
		\multicolumn{1}{r}{\textit{in}: 1} &\multicolumn{1}{r}{\textit{in}: 1}  & \multirow{2}{*}{2} & \multirow{2}{*}{14} \\
		\multicolumn{1}{r}{\textit{out}: 2} & \multicolumn{1}{r}{\textit{out}: 2} & &
	\end{tabular}
\end{table}

\subsection{3D-1D coupled problem with crossing inclusion} \label{test2}
The second numerical example takes into account a 3D-1D coupled problem, and proposes a validation of the XFEM strategy via a comparison with a solution obtained solving with standard FEM the original equi-dimensional problem. 

\begin{figure}
	\centering
	\begin{subfigure}[t]{0.45\textwidth}
		\centering
		\includegraphics[width=0.68\linewidth]{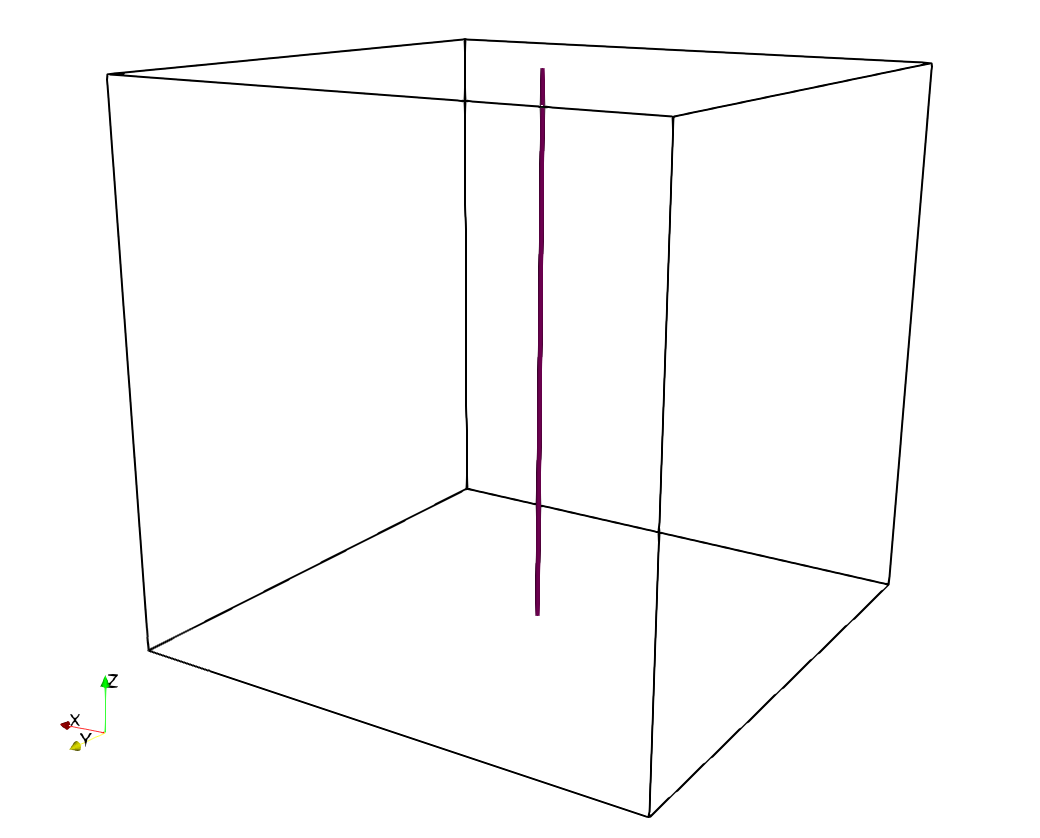}
		\caption{Geometry configuration}
		\label{1Dgeom_seg_passante}
	\end{subfigure}\hfill
	\begin{subfigure}[t]{0.45\textwidth}
		\hspace{0.5cm}
		\includegraphics[width=0.98\linewidth]{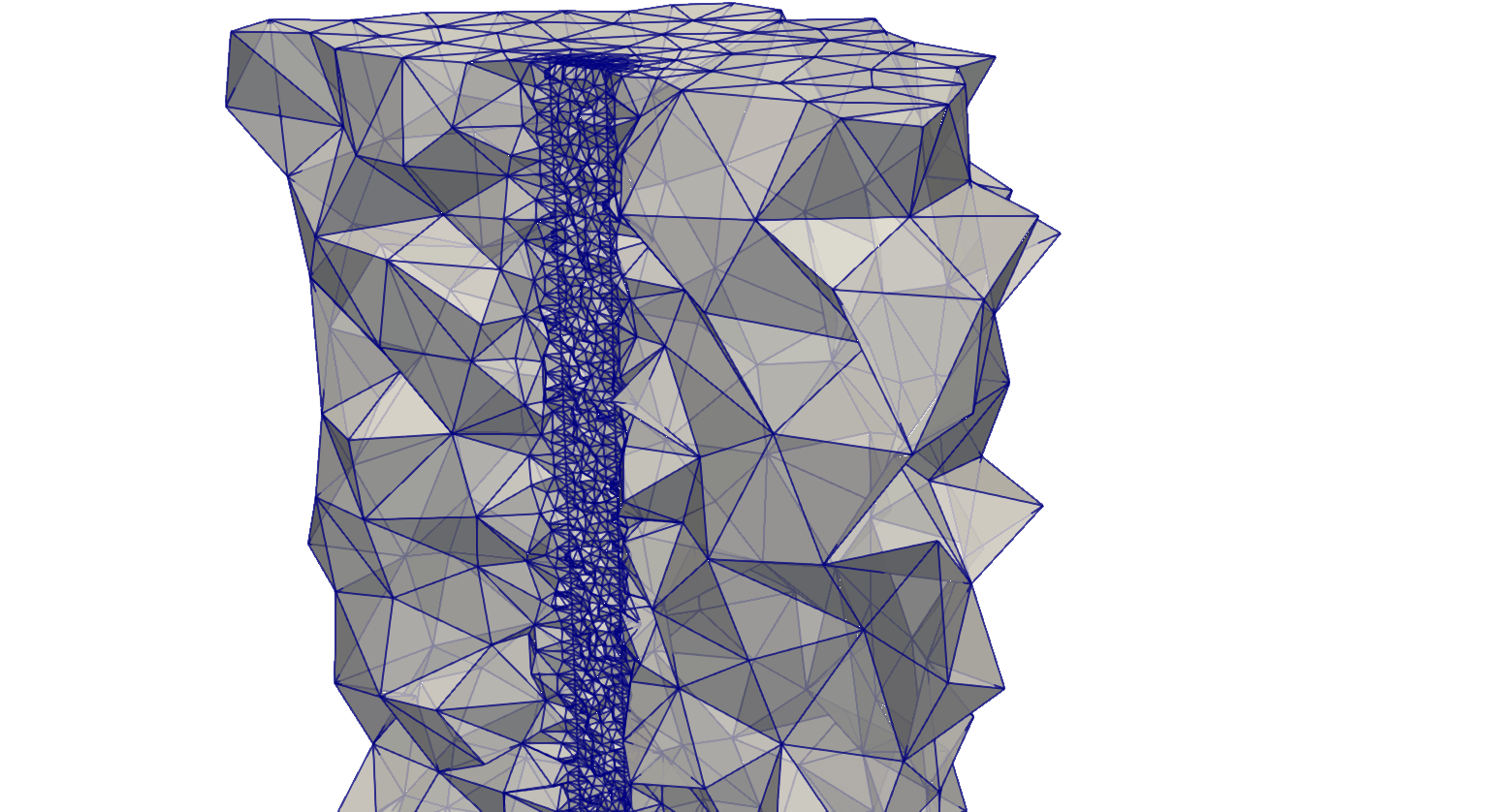}
		\caption{Mesh for the 3D-3D reference solution (detail)}
		\label{fig:3D3Dmesh_t2}
	\end{subfigure}
	\label{}
	\caption{Test \ref{test2}: Geometry configuration and detail of the mesh used for the reference solution.}
\end{figure}

Let us consider a cubic domain $\Omega=(-1,1)^3$ with a cylindrical inclusion $$\Sigma=\lbrace(x,y,z): \sqrt{x^2+y^2}<R,~ z\in (-1,1)\rbrace,$$ of radius $R=10^{-2}$ (see Figure \ref{1Dgeom_seg_passante}). We set $\partial \Omega_{\mathrm{d}}=\lbrace (x,y,z): z=-1 \vee z=1\rbrace$ and $\partial \Omega_{\mathrm{n}}=\partial \Omega \setminus \overline{\partial \Omega_{\mathrm{d}}}$, i.e., the Dirichlet boundary consists of the top and bottom faces of the cube, and the Neumann boundary consists of the lateral faces.
Problem data are $f=1$, $g=\overline{\overline{g}}=0$, $K=1$ and $\tilde{K}=10^5$, whereas homogeneous Dirichlet boundary conditions are prescribed on $\partial\Omega_{\mathrm{d}}$ and at the endpoints of $\Lambda$, while homogeneous Neumann boundary conditions are set on $\partial \Omega_{\mathrm{n}}$. 

To build a reference solution, we solve the original 3D-3D equi dimensional problem with a standard FEM method, on a mesh conforming to the actual interface $\Gamma$, which is discretized as the lateral surface of a prism with a 24-edge polygonal base. By standard FEM we actually mean that no domain decomposition is performed, and that a global pressure field $U$ is computed without resorting to an optimization based approach. As shown in Figure \ref{fig:3D3Dmesh_t2}, the mesh for the reference solution is refined in a region at a distance $R$ from $\Lambda$, with elements of maximum diameter of $0.002$, while it is coarser outside, where the element maximum diameter is $0.027$, resulting in about $3.1\times 10^5$ DOFs. 

The corresponding 3D-1D dimensionally reduced problem is solved on a uniformly refined mesh with element maximum diameter of $0.136$. We consider the cases $\rho=0$ and $\rho=0.01$, corresponding respectively to $N \sim 1.3\times10^3$ and $N\sim 1.5\times10^3$. Let us recall that the 3D-1D problem is always solved resorting to the optimization based domain decomposition method described in Section~\ref{sec:opt_form} and that, for $\rho>0$, we use the quadrature strategy described in Section \ref{quad} with the parameters reported in Table \ref{tab:quad_node}.

\begin{figure}
	\centering
	\includegraphics[width=0.4\linewidth]{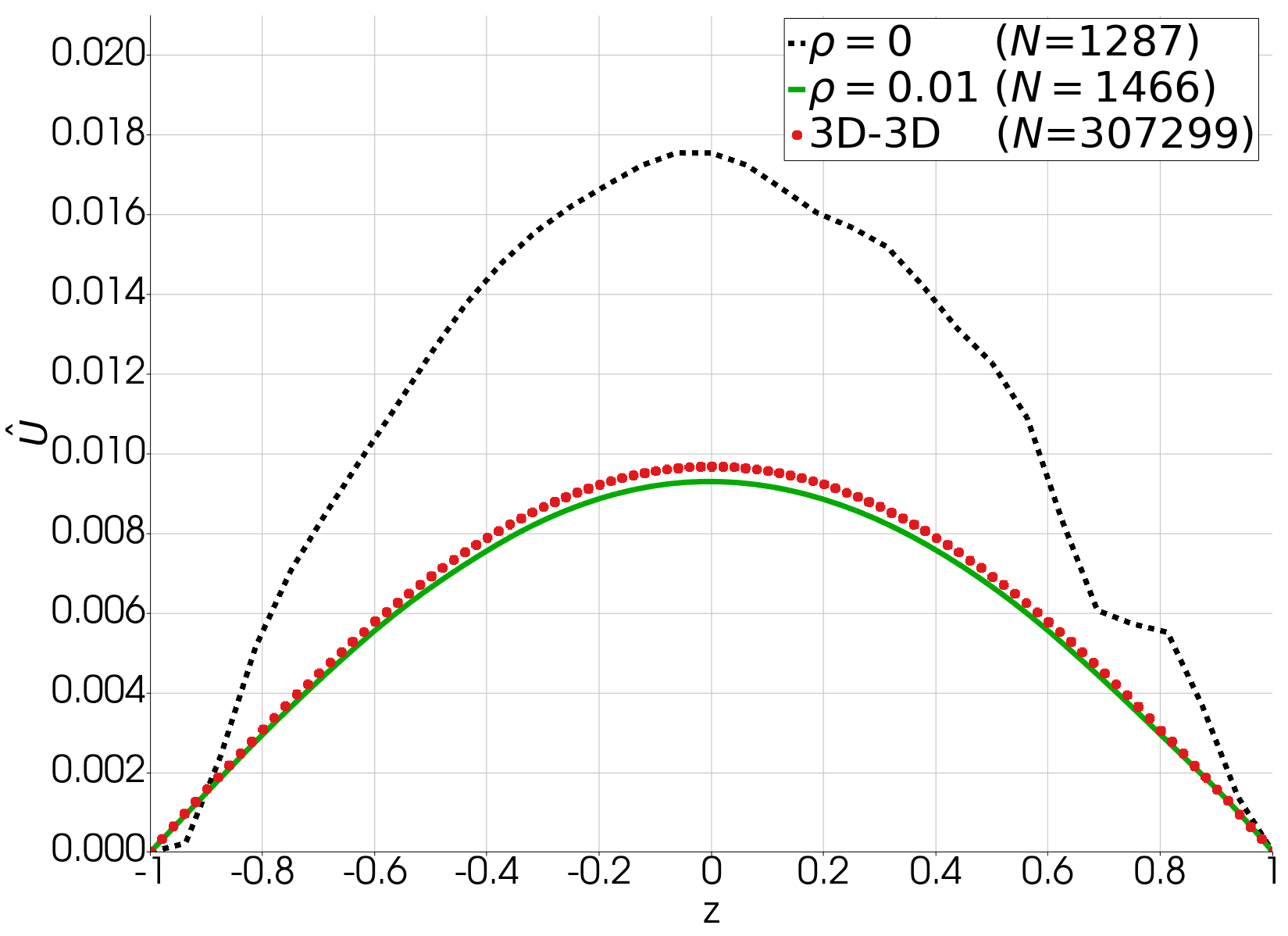}
	\caption{Test \ref{test2}: solutions on $\Lambda$ obtained for $\rho=0$ and $\rho=0.1$ compared to the trace on $\Lambda$ of the 3D-3D reference solution.}
	\label{fig:TC}
\end{figure}

The solutions obtained on $\Lambda$ are reported in Figure~\ref{fig:TC}, along with the trace on $\Lambda$ of the reference 3D-3D solution. In \cite{BGS3D1D2022}, where the problem was solved only for $\rho=0$, it was observed that, when $\tilde{K}\gg {K}$, mesh adaptation is needed close to the inclusion to improve accuracy, since a big jump in the diffusion coefficients produces a 3D solution with a very strong gradient close to the inclusion itself. Here we can instead see that, a choice of $\rho>0$ allows to obtain accurate solutions on uniform coarse meshes. Indeed, despite using nearly the same number of DOFs of $\rho=0$ case, the solution with the XFEM is almost overlapped to the 3D-3D reference solution.

\subsection{3D-1D coupled problem with inclusion inside the domain}
\label{test3}

\begin{figure}
	\centering
	\begin{subfigure}{.45\textwidth}
		\centering
		\includegraphics[width=0.89\linewidth]{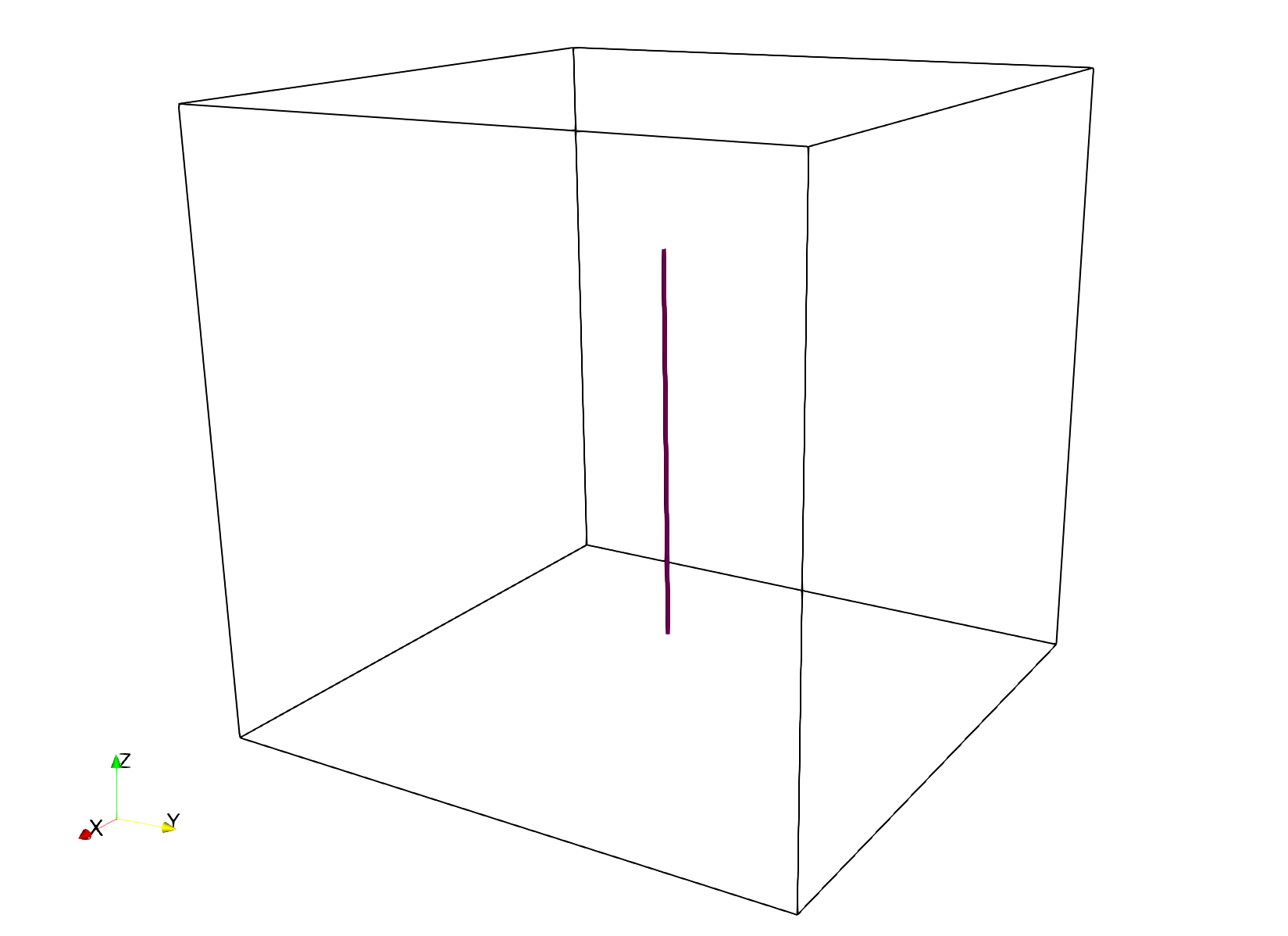}
		\caption{Geometry configuration}
		\label{fig:geom_t3}
	\end{subfigure}\hfill
	\begin{subfigure}{.45\textwidth}
		\centering
		\includegraphics[width=0.89\linewidth]{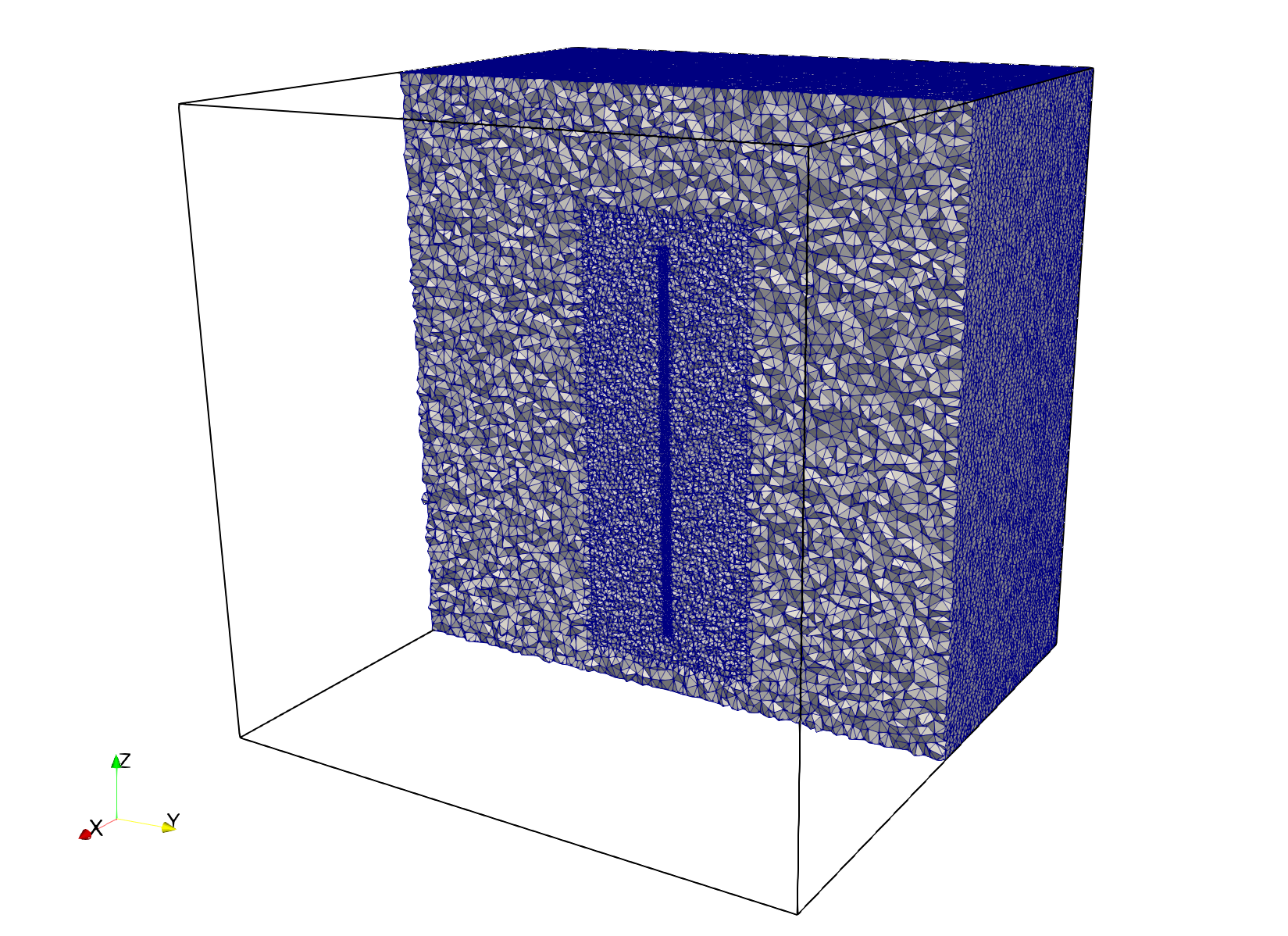}
		\caption{Mesh for the 3D-3D reference solution}
		\label{fig:adapted_mesh_t2}
	\end{subfigure}	
	\caption{Test~\ref{test3}: geometry configuration and mesh used to compute the reference solution}
	\label{fig:1dgeometry}
\end{figure}

In this example, we consider the case of an inclusion that is completely embedded into a domain $\Omega$. 
In particular we choose $\Omega=(-1,1)^3$ and the original 3D inclusion is
$$\Sigma=\lbrace(x,y,z): \sqrt{x^2+y^2}<R,~ z\in (-0.8,0.5)\rbrace,$$
i.e., the endpoints of $\Lambda$ lie inside $\Omega$, as reported in Figure \ref{fig:geom_t3}.
Problem data are $R=10^{-2}$, $f=1$, $K=1$, $g=\overline{\overline{g}}=0$, $\tilde{K}=10^5$ and we impose homogeneous Dirichlet boundary conditions on $\partial \Omega_{\mathrm{d}}=\lbrace (x,y,z): z=-1 \vee z=1\rbrace$ and homogeneous Neumann on $\partial \Omega_{\mathrm{n}}=\partial \Omega \setminus \overline{\partial \Omega_{\mathrm{d}}}$ and at the end sections of the inclusion.

As in the previous case, we build a reference solution by solving an equi-dimensional 3D-3D problem with standard finite elements on a mesh conforming to the interface $\Gamma$ and refined towards the inclusion
As reported in Figure \ref{fig:adapted_mesh_t2}, this mesh presents three different degrees of refinement. In particular the element maximum diameter is 0.0045 in a region of radius R around $\Lambda$, 0.01 outside this region but within a prismatic box $(-0.3,0.3)^2\times(-0.9,0.6)$, and 0.0215 outside the box, resulting in about $4 \cdot 10^5$ DOFs. 
The 3D-1D reduced problem in instead solved on a uniform mesh, with mesh parameter $0.086$, and for $\rho \in \{0,0.1,0.3,0.5,\sqrt{2}\}$, corresponding to a value of $N$ ranging between $4.6\times10^3$ and $1\times 10^4$.

\begin{figure}
	\centering
	\begin{subfigure}[t]{.48\textwidth}
		\centering
		\includegraphics[width=0.85\linewidth]{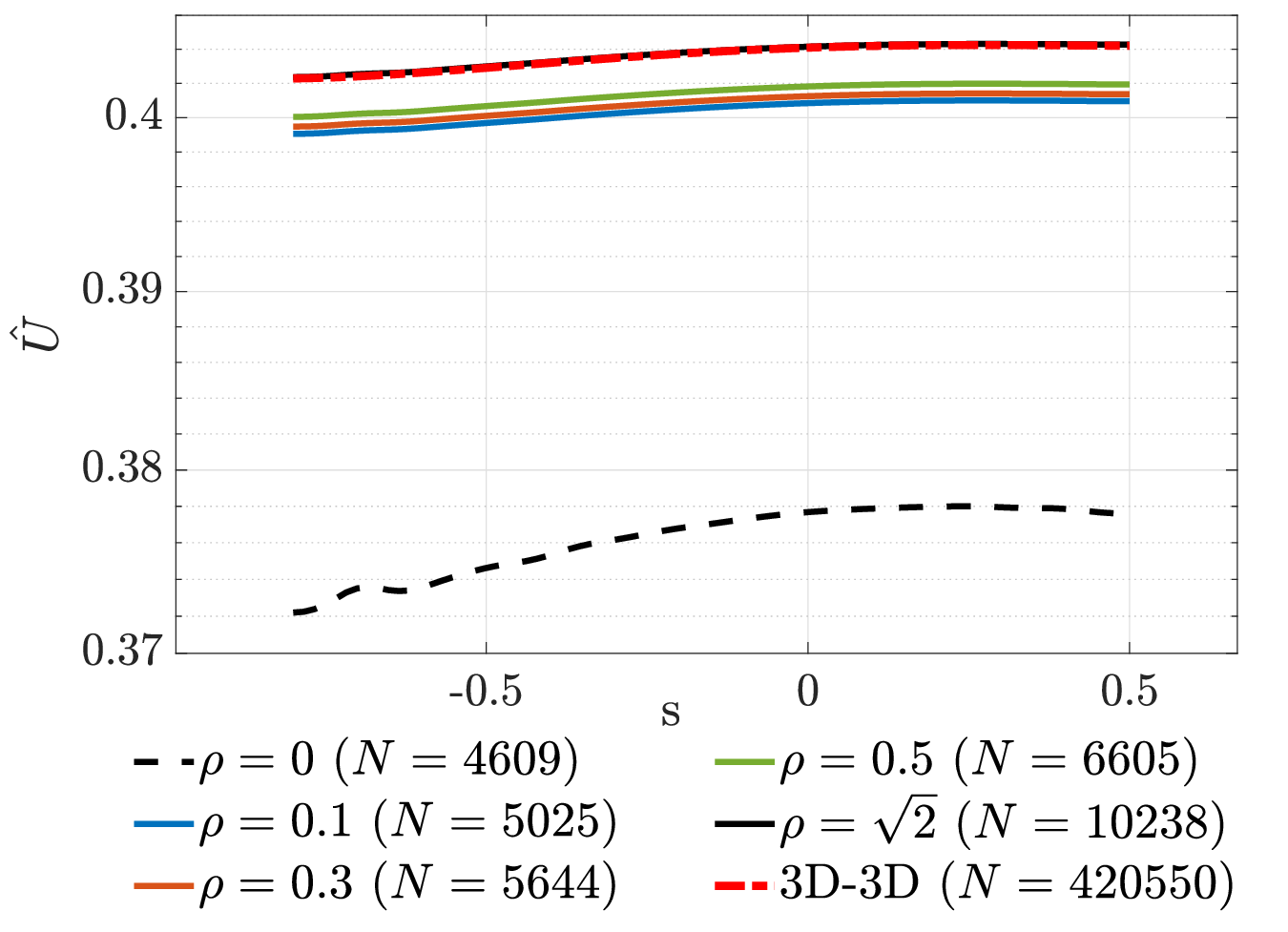}
		\caption{Solutions on $\Lambda$}
		\label{fig:internalendpoints-a}
	\end{subfigure}\hspace{0.5cm}%
	\hfill
	\begin{subfigure}[t]{.48\textwidth}
		\centering
		\includegraphics[width=0.85\linewidth]{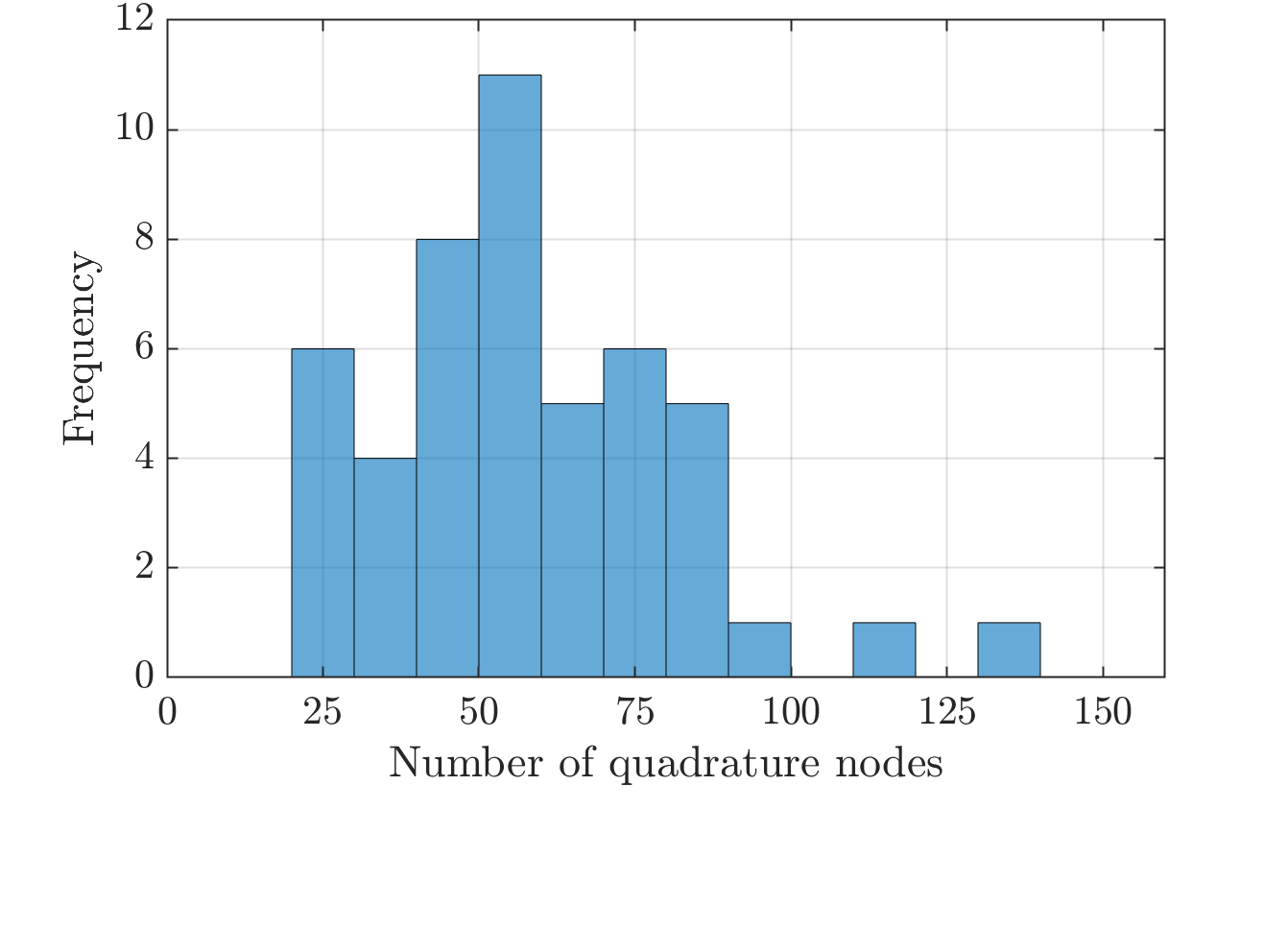}
		\caption{Distribution of quadrature nodes on elements intersected by $\Sigma$ when $\rho>0$.}
		\label{fig:internalendpoints-b}
	\end{subfigure}
	\caption{Test~\ref{test3}: solutions on $\Lambda$ and distribution of the number of quadrature nodes among elements intersected by the bulk inclusion $\Sigma$ when $\rho>0$. }
	\label{fig:internalendpoints}
\end{figure}

Figure~\ref{fig:internalendpoints-a} shows the solutions obtained on $\Lambda$. We can observe that for $\rho>0$, the solutions are much closer to the trace of the reference 3D-3D solution, already for $\rho=0.1$. Clearly, as $\rho$ increases, the gap with the reference decreases, at the expenses of a larger number of unknowns. The case $\rho=0$, instead, fails in providing a good representation of the solution, since, as in the previous test, the large jump in the coefficients between the 3D domain and the 1D inclusion gives a solution with a steep gradient that can not be correctly reproduced by FEM basis functions on elements with a diameter larger than the radius of the inclusion. We remark that it is not possible to significantly reduce the number of DOFs of the 3D-3D reference solution without affecting its quality. We can then note that choosing $\rho=\sqrt{2}$, i.e. enriching all the basis functions, still gives a number of unknowns about $40$ times smaller than the ones required for the equi-dimensional problem.

Figure~\ref{fig:internalendpoints-b} reports the distribution of the total number of quadrature nodes used in the elements intersected by the inclusion when $\rho>0$. The maximum number of quadrature nodes in a non-split cell can be easily computed as $3\times 4\times n_\Lambda \times n_\mathfrak{r} \times n_\theta=100$, since $3$ is the maximum number of intervals $\Lambda^t$ in a tetrahedron, $4$ is the maximum number of triangular regions on each slice of Figure~\ref{Fig_quad} and $n_\Lambda$, $n_\mathfrak{r}$, and $n_\theta$ are the values taken from Table~\ref{tab:quad_node}, summing the nodes inside and outside $\Gamma$. This is confirmed by the values in Figure~\ref{fig:internalendpoints-b}, with the only exception of the two elements containing inclusion endpoints, that are, instead, split into sub-cells. Clearly the high number of quadrature nodes represents an additional computational cost. However, in general, this  cost is largely offset by the possibility of using less degrees of freedom with respect to approaches that require mesh adaptation. Moreover, the quadrature rule described in Section~\ref{quad} is only used in elements intersected by $\Sigma$, and thus the values of Figure~\ref{fig:internalendpoints-b} are independent of the chosen value of $\rho>0$. When the mesh-size is reduced, the number of such elements grows linearly as $h^{-1}$.

\subsection{3D-1D coupled problem with bifurcated inclusion}\label{test4}
\begin{figure}
	\begin{subfigure}[t]{.45\textwidth}
		\centering
		\includegraphics[width=0.85\linewidth]{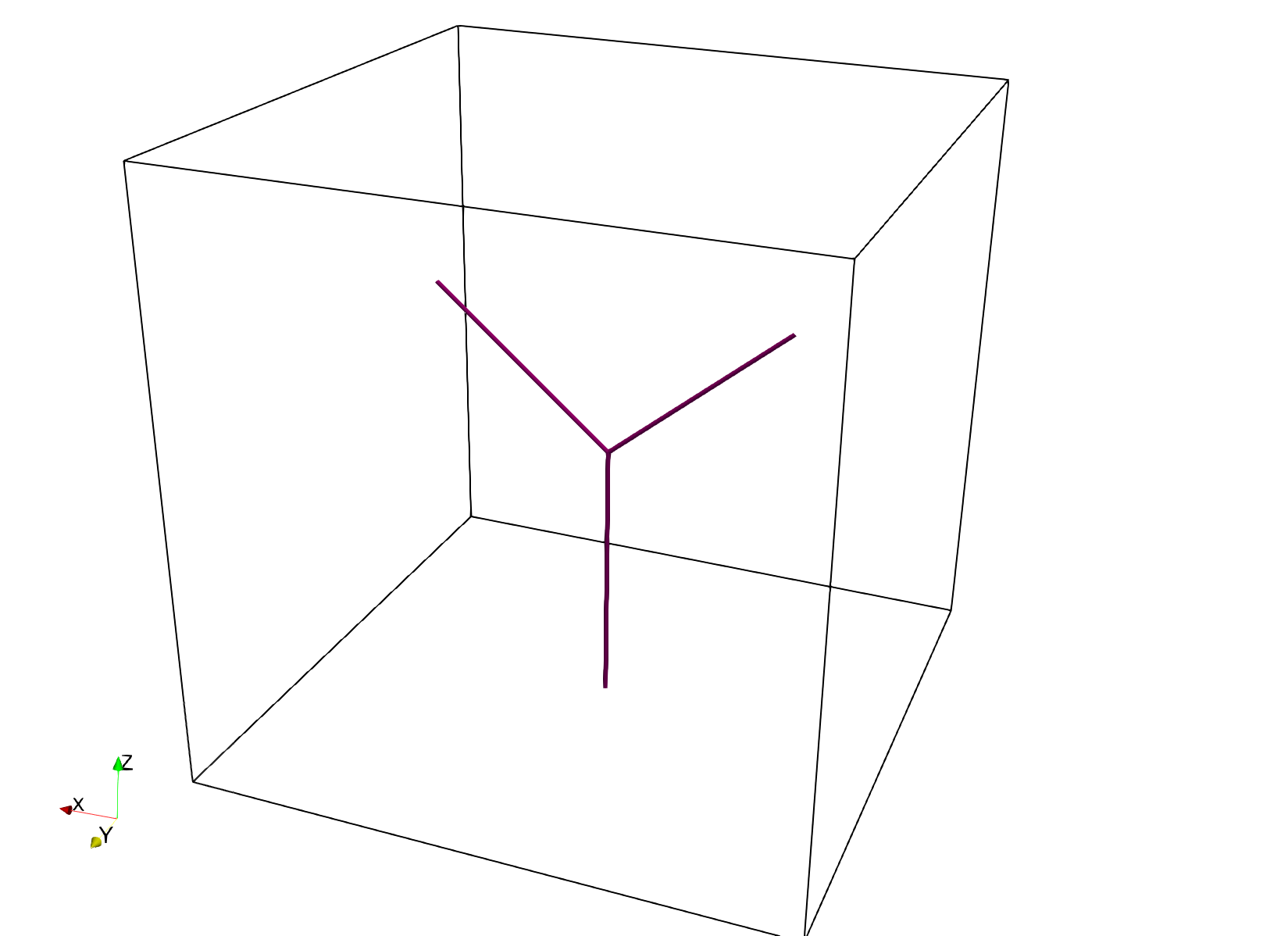}
		\caption{Geometry configuration}
		\label{fig:1Dgeom}
	\end{subfigure}\hfill
	\begin{subfigure}[t]{.45\textwidth}
		\centering
		\includegraphics[width=0.85\linewidth]{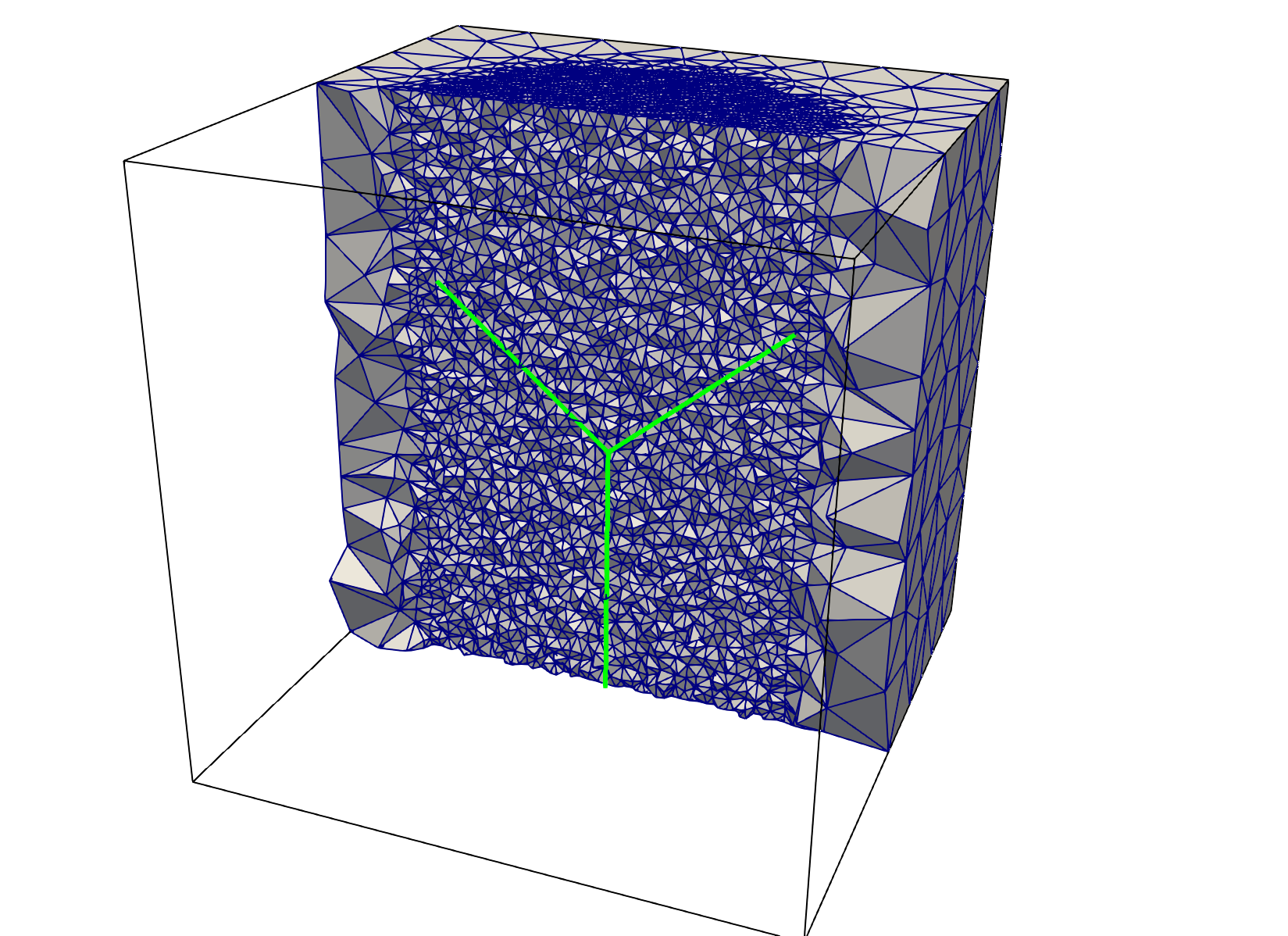}
		\caption{Mesh for the reference solution}
		\label{fig:test4_mesh_adapted}
	\end{subfigure}\hfill
	\caption{Test \ref{test4}: geometry configuration and mesh used to compute the reference solution.}
	\label{}
\end{figure}

Let us now consider the case of a bifurcated inclusion $\Sigma$, which can also be seen as the case of multiple inclusions $\Sigma_i$ whose centerlines $\Lambda_i$ intersect at one point.
In particular, we consider $3$ inclusions of radius $R=10^{-2}$ with centrelines
\begin{align*}
&\Lambda_1=\lbrace(0,0,z): z \in (-1.0,-0.1)\rbrace\\
&\Lambda_2=\lbrace(x,0,z): x \in (0,0.6),~z \in (-0.1,0.4)\rbrace\\
&\Lambda_3=\lbrace(x,0,z): x \in (0,-0.6),~z \in (-0.1,0.4)\rbrace,
\end{align*}
as shown in Figure~\ref{fig:1Dgeom}.
Let $\Omega=(-1,1)^3$ and let us enforce homogeneous Dirichlet boundary conditions on the top and bottom faces of the cube. 
We also impose a homogeneous Dirichlet boundary condition on the section of the inclusion lying on the bottom face of $\Omega$, and homogeneous Neumann boundary conditions on the sections lying inside $\Omega$. We finally set $f=1$, $K=1$, $g=\overline{\overline{g}}=0$ and $\tilde{K}=10^5$. 

To obtain a reference solution, in this case, we solve the 3D-1D reduced problem with $\rho=0$, but on a mesh refined within a prism containing the whole inclusion (see Figure \ref{fig:test4_mesh_adapted}). In particular, the prism has a 7-edge polygonal base which can be inscribed in a circle of radius 0.7. Inside the prism the mesh parameter is $0.027$, while it is $0.215$ outside, resulting in about $4.6\times 10^4$ DOFs. This choice of reference solution is not as reliable as the one of the previous examples, and is motivated by the complexity of generating a mesh conforming to the 3D inclusion for complex geometries, as the ones proposed here and in example \ref{test_final}.

The reduced 3D-1D problem is then solved on a uniform mesh with mesh parameter $0.086$ and for $\rho \in \{0,0.1, 0.3,0.5,\sqrt{2}\}$, corresponding to $N \in[4.6\times 10^3, 2.2 \times 10^4]$. Let us remark that we are choosing a unique value of $\rho$ for all the inclusions. In the following, when specifying the value of $\rho$, we will always refer to a solution computed on the uniform mesh.

\begin{figure}
	\centering
	\begin{subfigure}{.33\textwidth}
		\centering
		\includegraphics[width=1\linewidth]{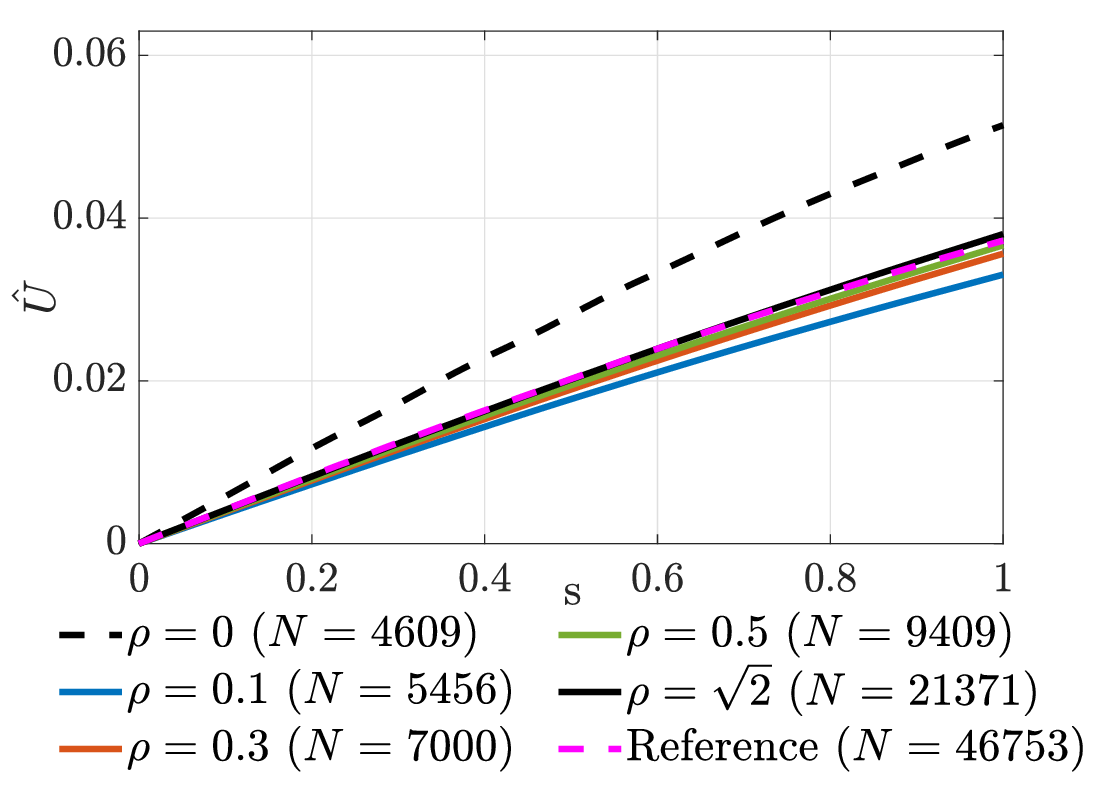}
		\caption{$\Lambda_1$}
		\label{fig:1dsolbranch-a}
	\end{subfigure}\hspace{0.1cm}%
	\begin{subfigure}{.33\textwidth}
		\centering
		\includegraphics[width=1\linewidth]{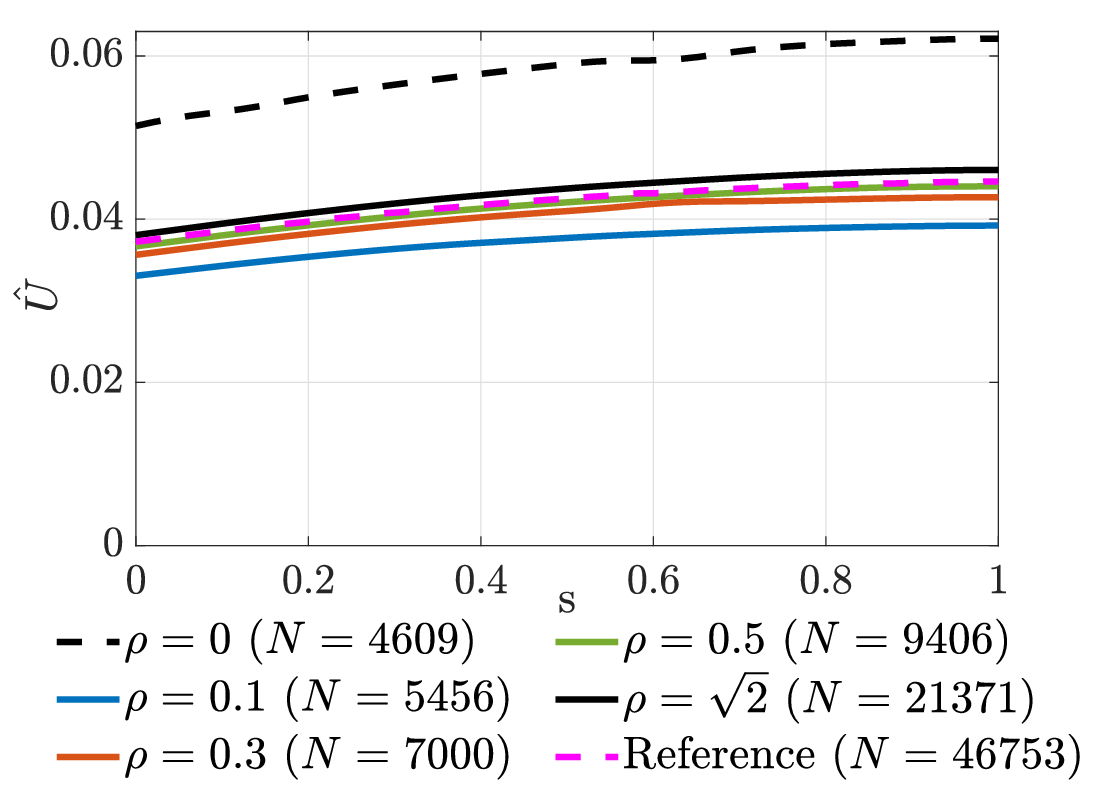}
		\caption{$\Lambda_2$}
		\label{fig:1dsolbranch-b}
	\end{subfigure}\hspace{0.05cm}%
	\begin{subfigure}{.33\textwidth}
		\centering
		\includegraphics[width=1\linewidth]{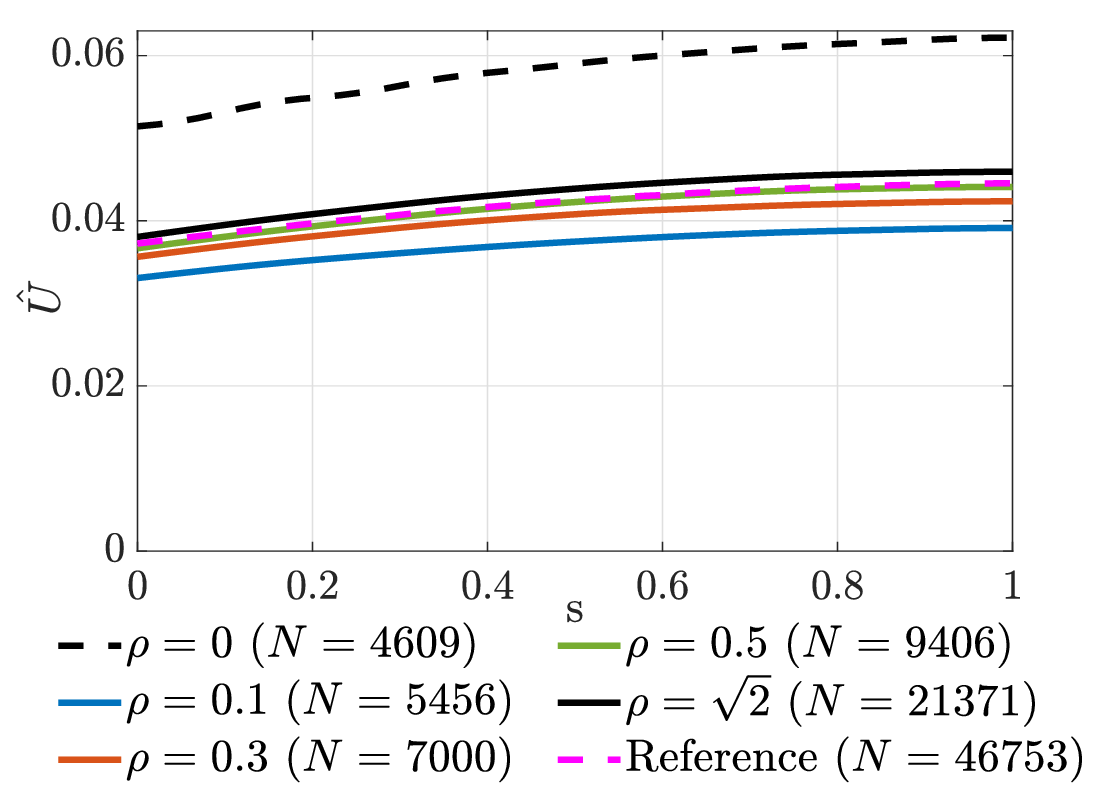}
		\caption{$\Lambda_3$}
		\label{fig:1dsolbranch-c}
	\end{subfigure}
	\caption{Test \ref{test4}: solutions on $\Lambda$ for different values of $\rho$ compared to reference solution.}
	\label{fig:1dsolbranch}
\end{figure}

Figure~\ref{fig:1dsolbranch} shows the solutions for the different values of $\rho$ on $\Lambda_1$, $\Lambda_2$, $\Lambda_3$, along with the corresponding trace of the reference solution. Also in this case we can notice that $\rho>0$ provides good approximations of the reference solution, whereas the case $\rho=0$ is less accurate on the uniform mesh. 

\begin{figure}
	\begin{subfigure}[t]{.47\textwidth}
		\centering
		\includegraphics[width=0.75\linewidth]{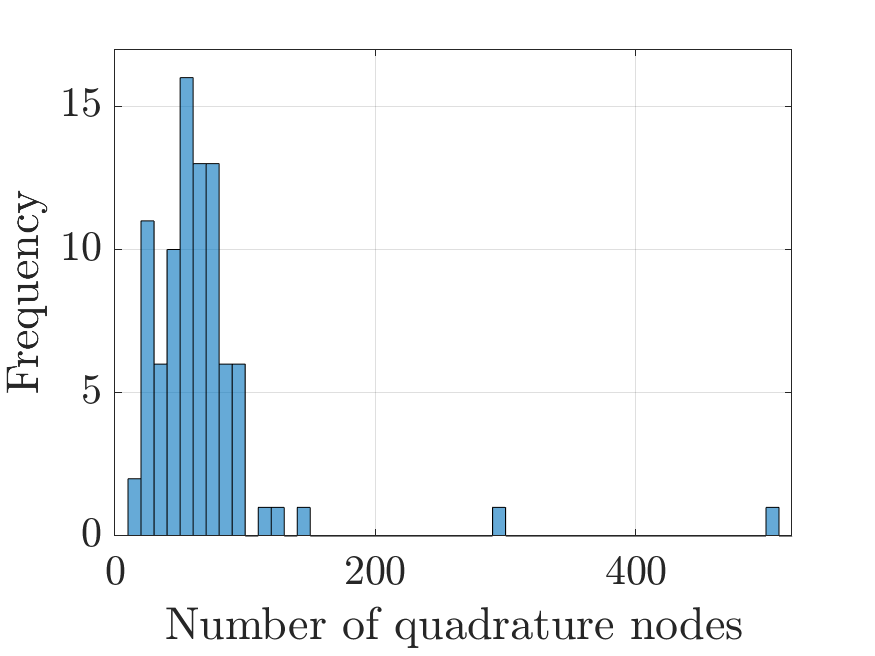}
		\caption{Distribution of quadrature nodes on elements intersected by $\Sigma$ for $\rho>0$.}
		\label{fig:qnodes}
	\end{subfigure}\hfill
	\begin{subfigure}[t]{.47\textwidth}
		\centering	\includegraphics[width=0.75\linewidth]{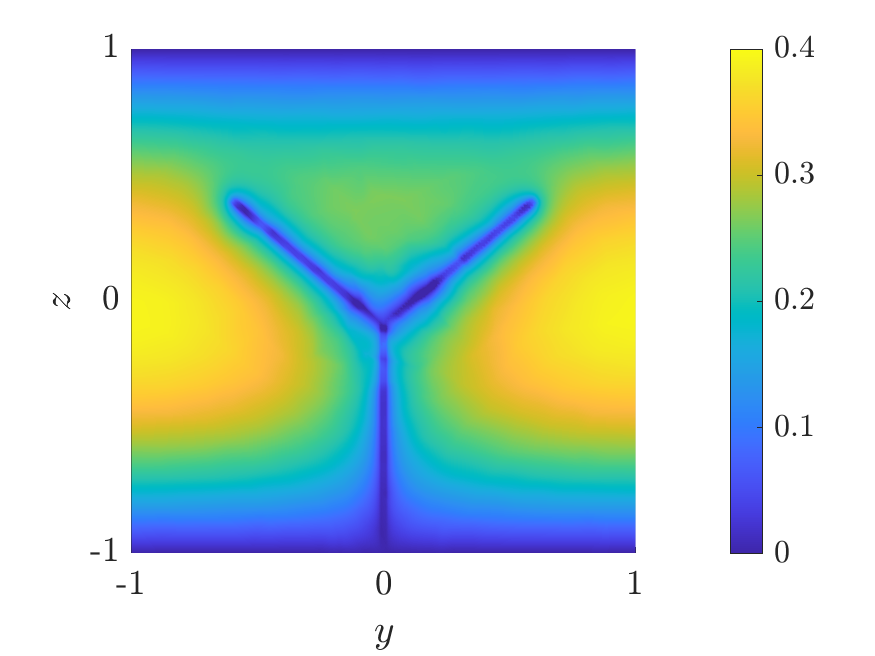}
		\caption{Solution on the $yz$-plane for $\rho=0.1$}
		\label{fig:3dsolbranch_yz}
	\end{subfigure}
	\caption{Test \ref{test4}: distribution of the number of quadrature nodes among elements intersected by the bulk inclusion $\Sigma$ when $\rho>0$ and solution on the $yz$-plane obtained for $\rho=0.1$.}
	\label{}
\end{figure}

The distribution of the total number of quadrature nodes used in the elements intersected by $\Sigma$ when $\rho>0$ is reported in Figure~\ref{fig:qnodes}. In this case only four elements exceed $100$ quadrature nodes: the two containing the endpoints of the inclusion, the one containing the intersection point and one of its neighbors, which are the elements which actually require splitting. Figure~\ref{fig:3dsolbranch_yz} shows a slice of the solution obtained for $\rho=0.1$ on the $yz$ plane.

\begin{figure}
	\centering
	\begin{subfigure}[t]{.5\textwidth}
		\centering
		\includegraphics[width=0.8\linewidth]{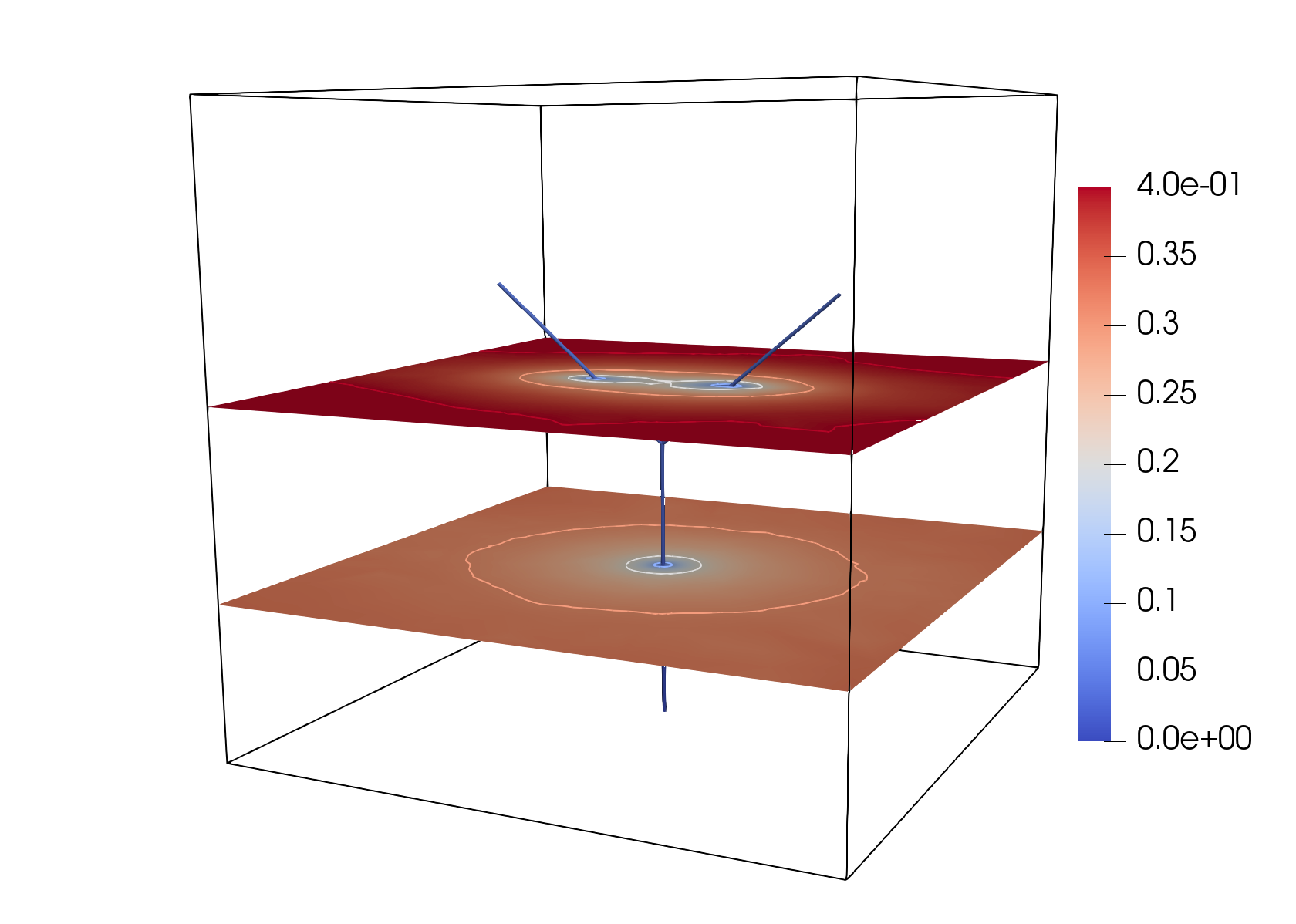}
		\caption{Solution for $\rho=0.1$ ($N=5456$).}
		\label{fig:3dsolbranch_slice-a}
	\end{subfigure}%
	\hfill
	\begin{subfigure}[t]{.5\textwidth}
		\centering
		\includegraphics[width=0.8\linewidth]{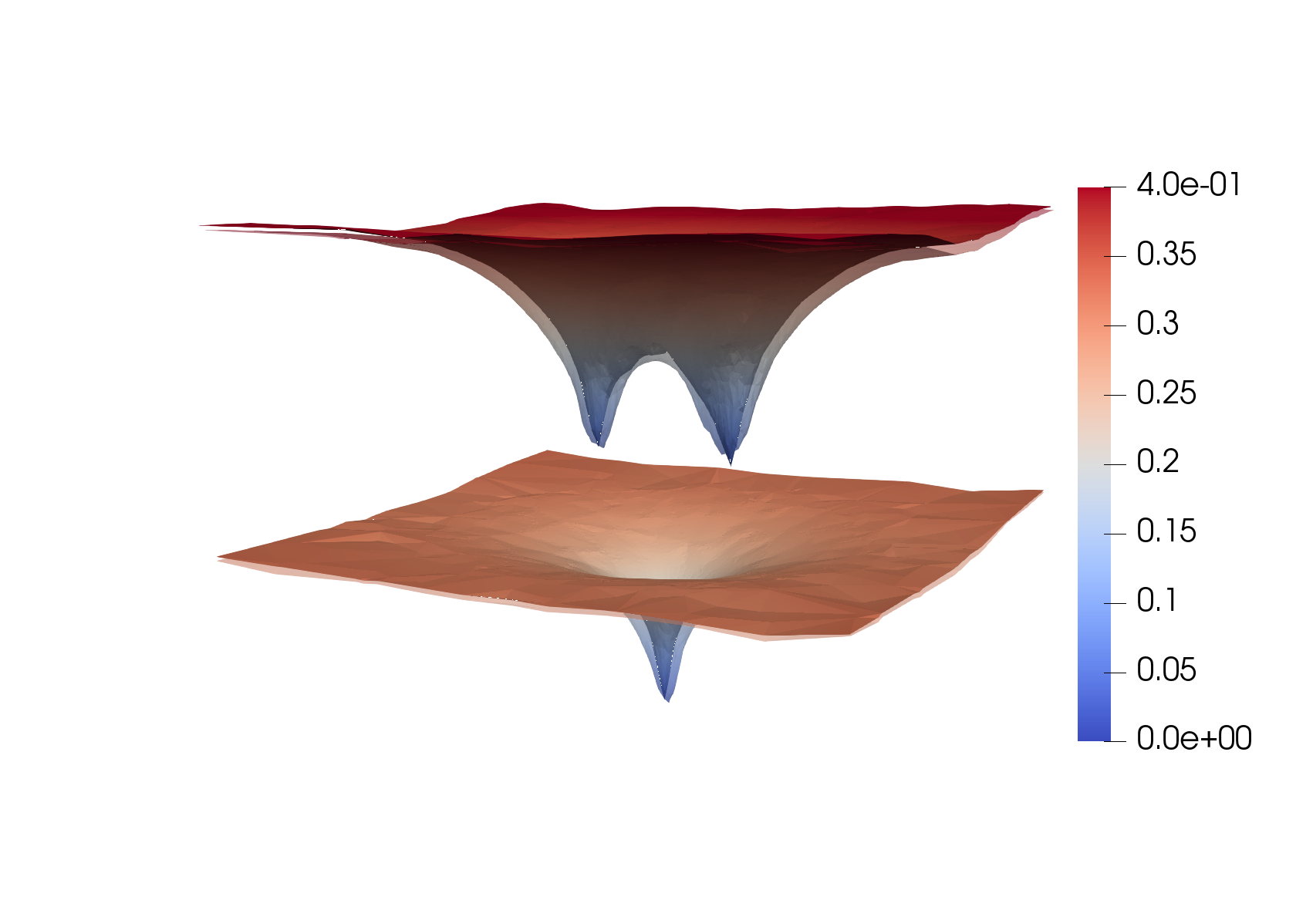}
		\caption{Solution obtained for $\rho=0.1$ ($N=5456$, opaque) and reference solution ($N=46753$, transparent)}
		\label{fig:3dsolbranch_slice-b}
	\end{subfigure}
	\caption{Test \ref{test4}: solution on planes orthogonal to the $z$-axis and located at $z=-0.5$ and $z=0.1$.}
	\label{fig:3dsolbranch_slice}
\end{figure}

\begin{figure}
	\centering
	\begin{subfigure}[t]{.5\textwidth}
		\centering
		\includegraphics[width=0.8\linewidth]{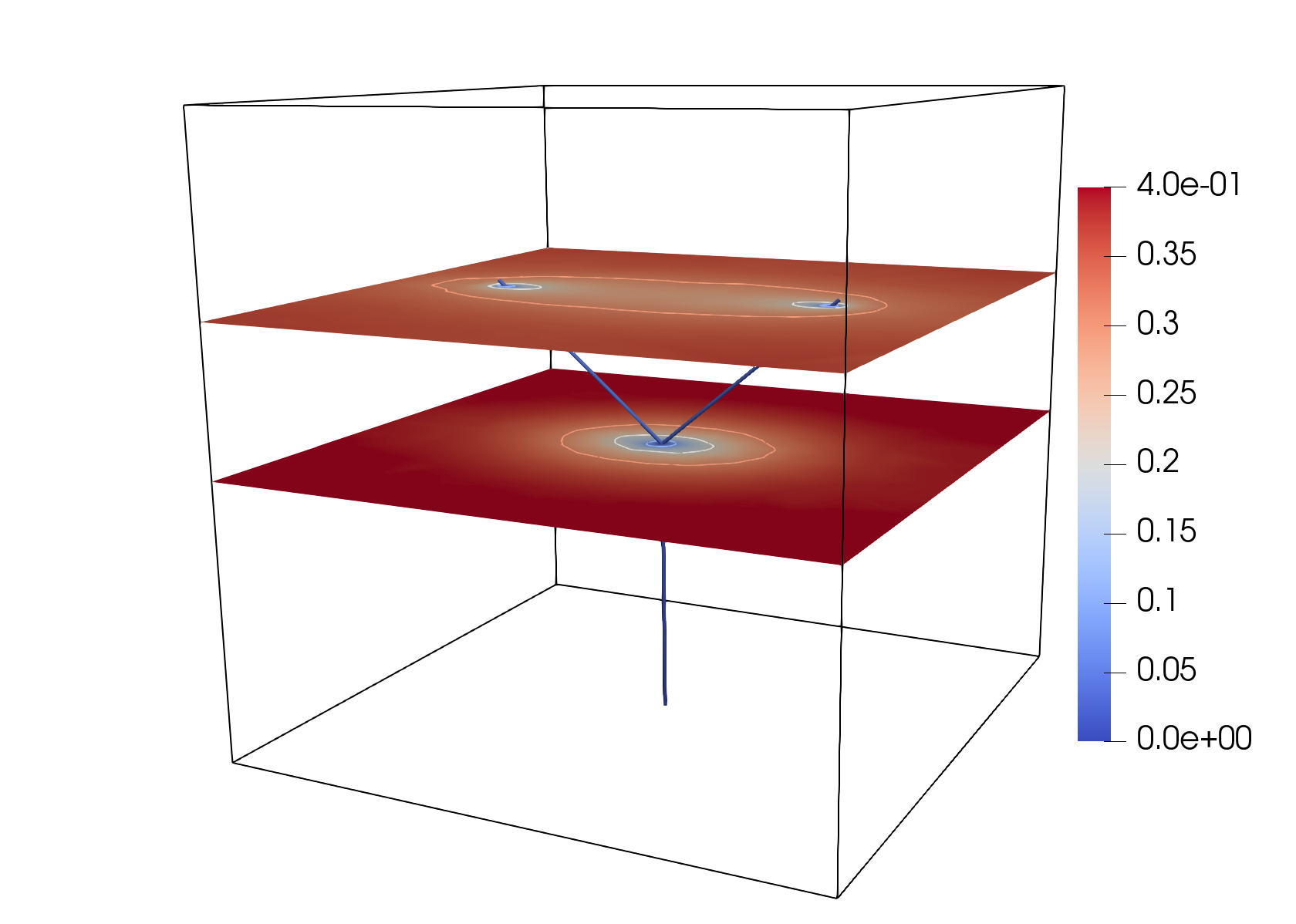}
		\caption{Solution for $\rho=0.1$ ($N=5456$).}
		\label{fig:3dsolbranch_slice_tips-a}
	\end{subfigure}%
	\hfill
	\begin{subfigure}[t]{.5\textwidth}
		\centering
		\includegraphics[width=0.8\linewidth]{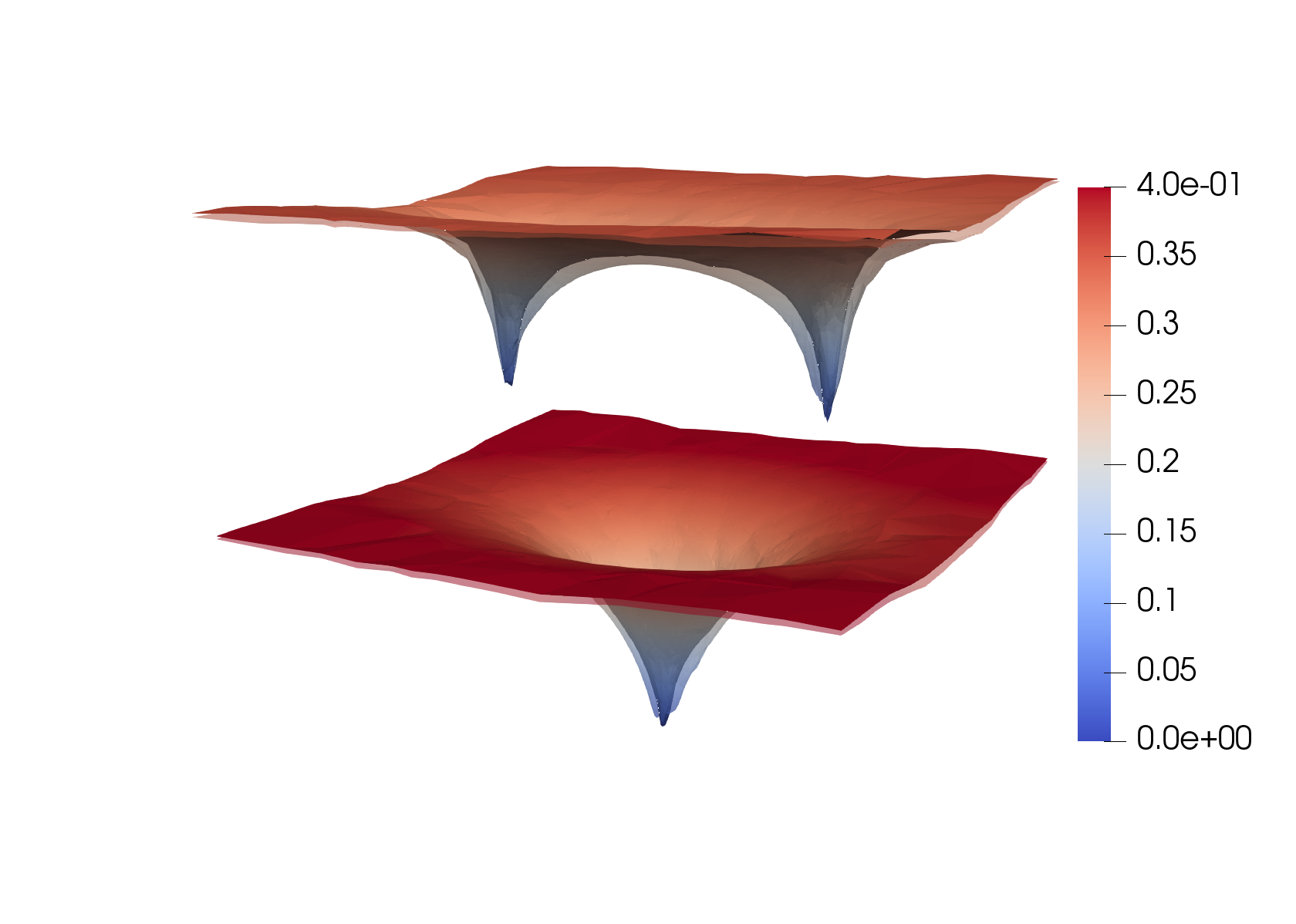}
		\caption{Solution obtained for $\rho=0.1$ ($N=5456$, opaque) and reference solution ($N=46753$, transparent)}
		\label{fig:3dsolbranch_slice_tips-b}
	\end{subfigure}
	\caption{Test \ref{test4}: solution on planes orthogonal to the $z$-axis and passing through the bifurcation point and close to the two tips.}
	\label{fig:3dsolbranch_slice_tips}
\end{figure}

Finally, Figures~\ref{fig:3dsolbranch_slice}-\ref{fig:3dsolbranch_slice_tips} show the solution for $\rho=0.1$ on planes orthogonal to the $z$-axis. In the right panels, the reference solution is reported in transparency, to highlight the good matching of the two.


\subsection{3D-1D coupled problems: inclusion with several branches}\label{test_final}
\begin{figure}
	\centering
	\begin{subfigure}[t]{.45\textwidth}
		\centering
		\includegraphics[width=0.85\linewidth]{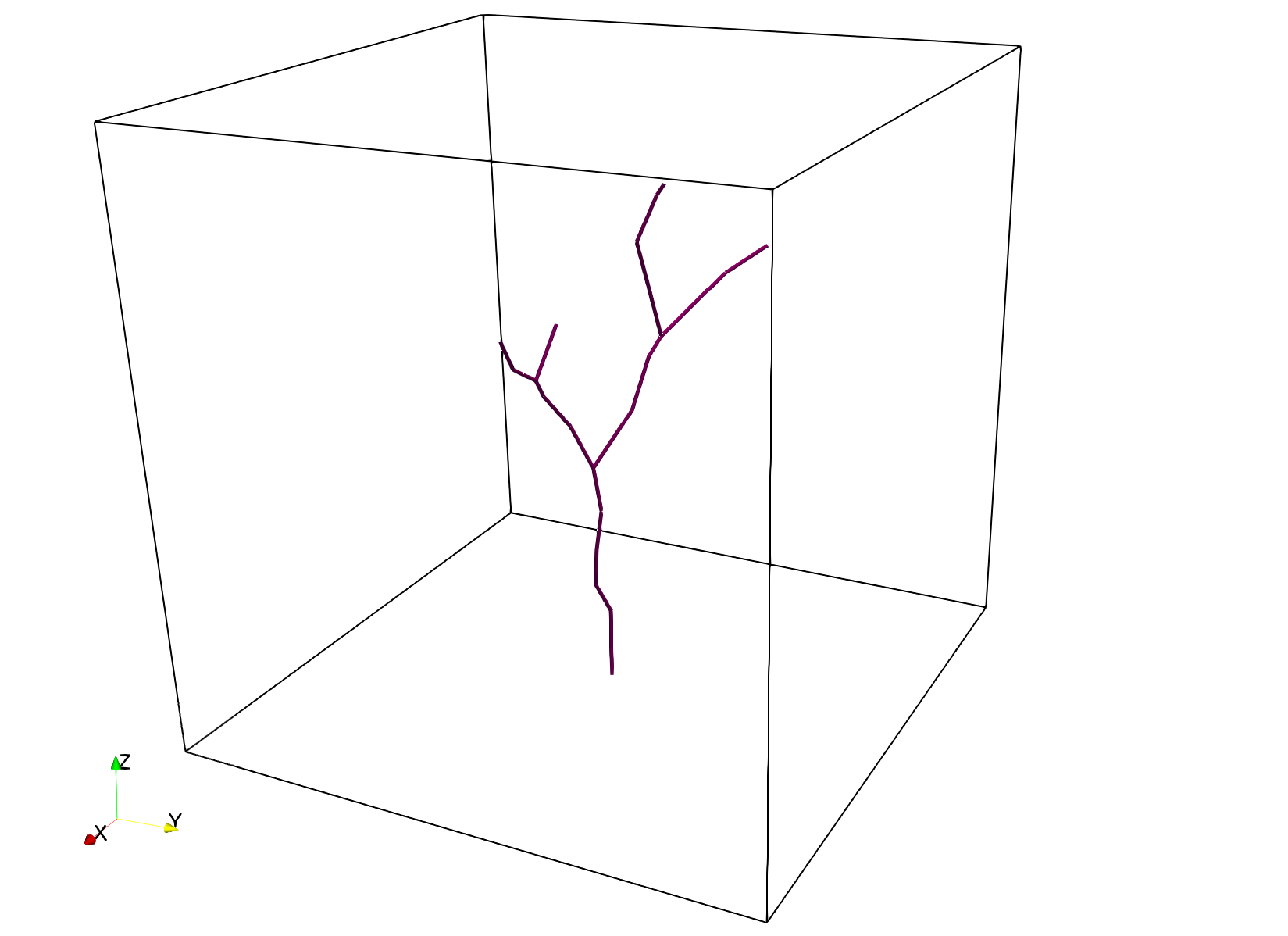}
		\caption{Geometry configuration}
		\label{fig:geom}
	\end{subfigure}\hfill 
	\centering
	\begin{subfigure}[t]{.45\textwidth}
		\centering
		\includegraphics[width=0.85\linewidth]{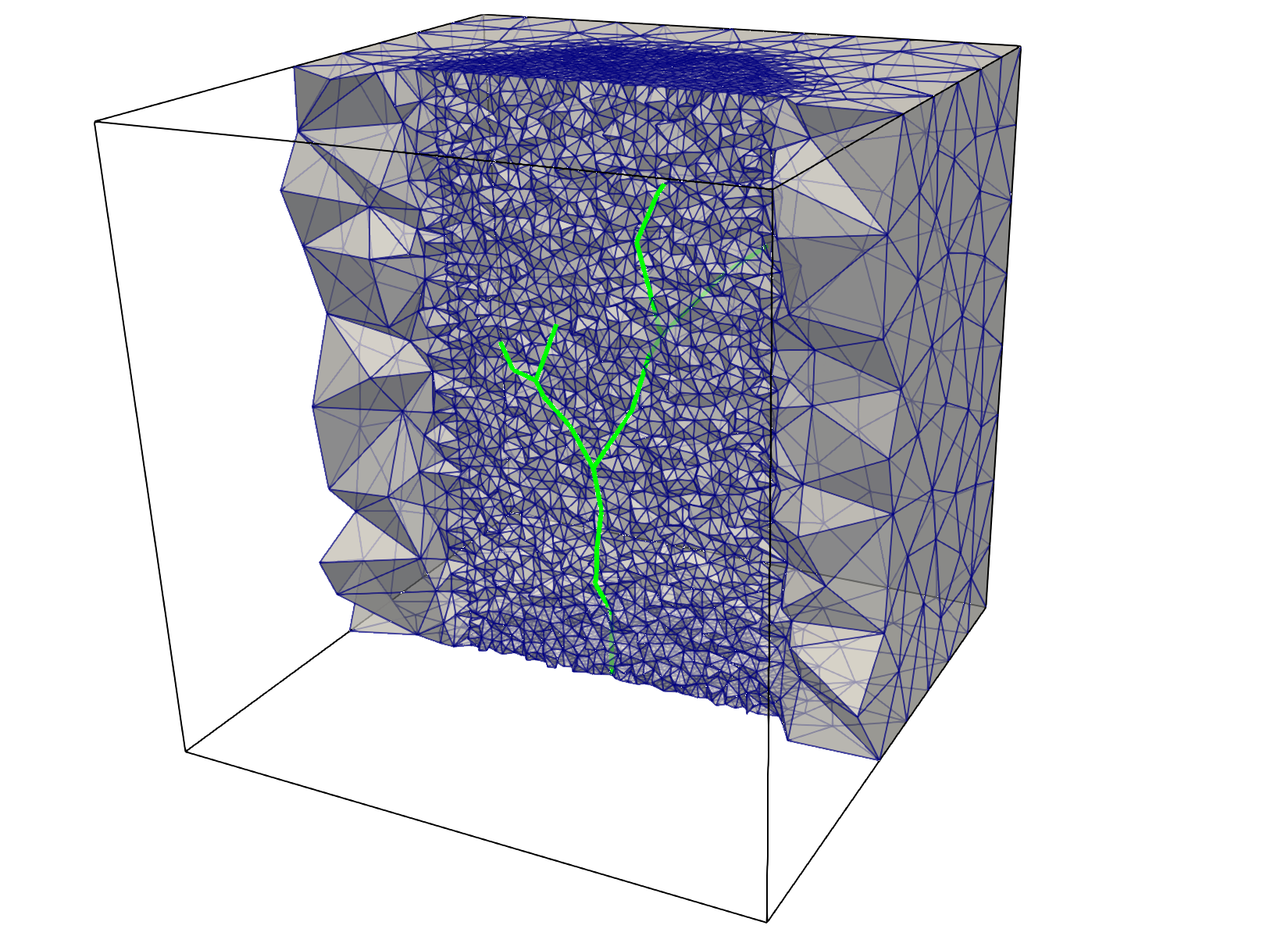}
		\caption{Refined mesh for the reference solution}
		\label{fig:test5_mesh_adapted}
	\end{subfigure}
	\caption{Test \ref{test_final}: geometry configuration and mesh used to compute the reference solution.}
	\label{fig:geom_and_mesh}
\end{figure}

As a last numerical example, we propose a case with a more realistic inclusion characterized by several branches, as reported in Figure \ref{fig:geom} . We assume that the inclusion $\Sigma$, which has a constant radius $R=10^{-2}$, is embedded in a cubic domain $\Omega=(-1,1)^3$. We chose $f=0$, $K=1$, $g=\overline{\overline{g}}=0$, $\tilde{K}=10^5$ and we impose homogeneous Neumann boundary conditions on $\partial \Omega_{\mathrm{n}}=\lbrace (x,y,z): z=-1 \vee z=1\rbrace$ and homogeneous Dirichlet on $\partial \Omega_{\mathrm{d}}=\partial \Omega \setminus \partial \Omega_{\mathrm{n}}$. For whats concerns the inclusion end sections, we prescribe a Dirichlet boundary condition equal to one at the section lying on the bottom face of the cube, while homogeneous Neumann conditions are prescribed at the dead ends.

As for the previous test case, we build a reference solution by solving the 3D-1D reduced problem with $\rho=0$ on a mesh refined in a prism containing the inclusion. In this case we consider a prism with a 7-edge polygonal base which can be inscribed in a circle of radius 0.6. Inside the prism we consider a mesh parameter of $0.027$, while outside of 0.215, resulting in $N\sim4.3\times 10^4$ (see Figure \ref{fig:test5_mesh_adapted}).

The 3D-1D reduced problem is solved on a uniform mesh of parameter $0.086$ with $\rho=0$ and $\rho=0.1$, resulting in $N\sim3.3\times 10^3$ and $N\sim4.6\times 10^3$ respectively. Also in this case, if the value of $\rho$ is specified, we will always refer to a solution computed on the uniform mesh.

\begin{figure}
	\centering
	\includegraphics[width=0.4\linewidth]{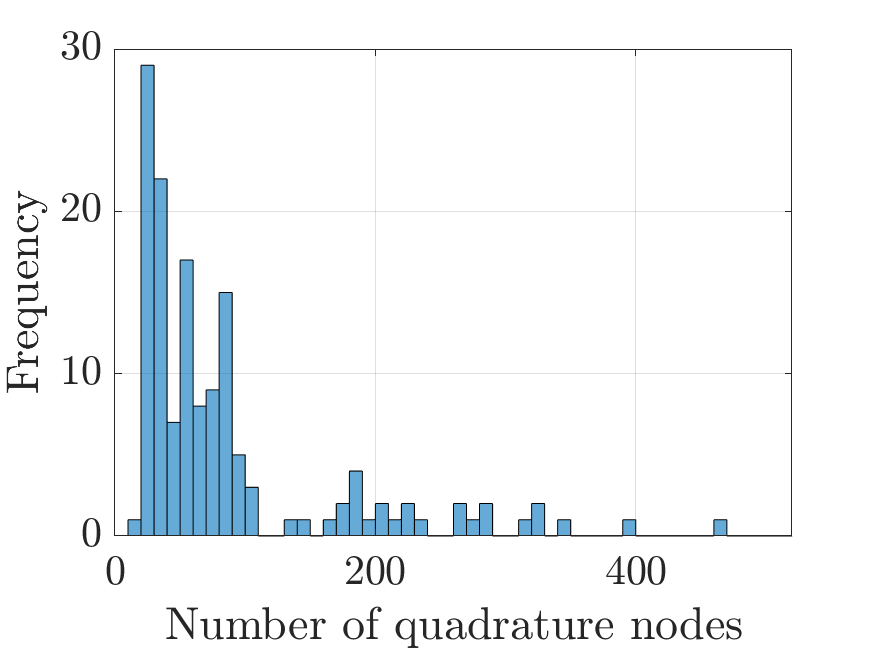}
	\caption{Test \ref{test_final}: distribution of the number of quadrature nodes among the elements intersected by the bulk inclusion $\Sigma$ when $\rho>0$.}
	\label{fig:test5_qnodes}
\end{figure}

The distribution of the number of quadrature nodes among the elements cut by $\Sigma$ when $\rho>0$ is reported in Figure \ref{fig:test5_qnodes}, using again the quadrature parameters reported in Table \ref{tab:quad_node}. As for the previous cases, only a few tetrahedrons present a particularly high number of quadrature nodes, i.e. the ones containing the bifurcation points or being close to them, and the ones containing the dead ends.

\begin{figure}
	\centering
	\begin{subfigure}{.45\textwidth}
		\centering
		\includegraphics[width=1.0\linewidth]{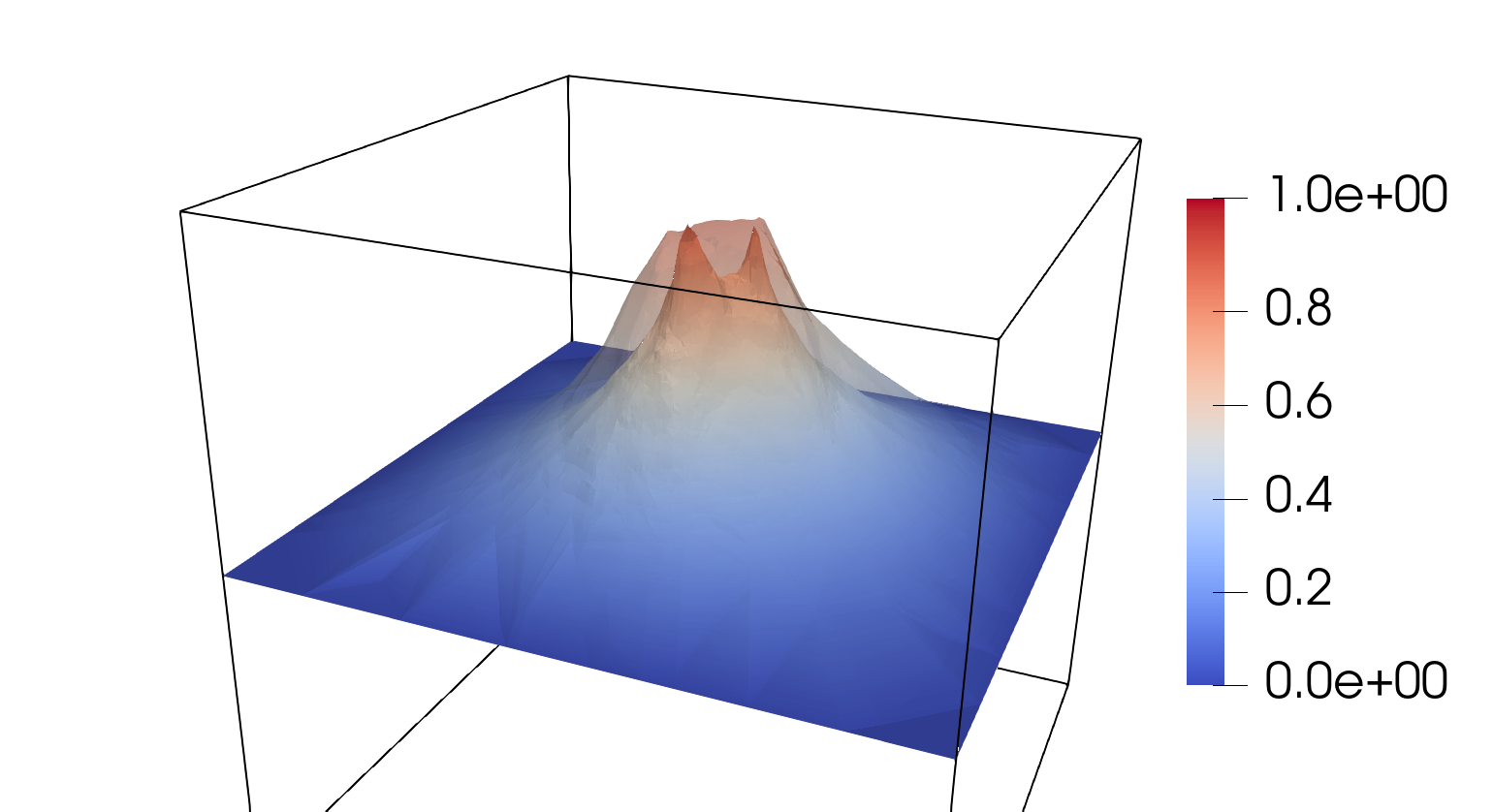}
		\caption{Solution obtained for $\rho=0$ (transparent, $N=3331$) and for $\rho=0.1$ (opaque, $N=4574$). }
		\label{fig:warp-a}
	\end{subfigure}\hfill 
	\centering
	\begin{subfigure}{.45\textwidth}
		\centering
		\includegraphics[width=1.0\linewidth]{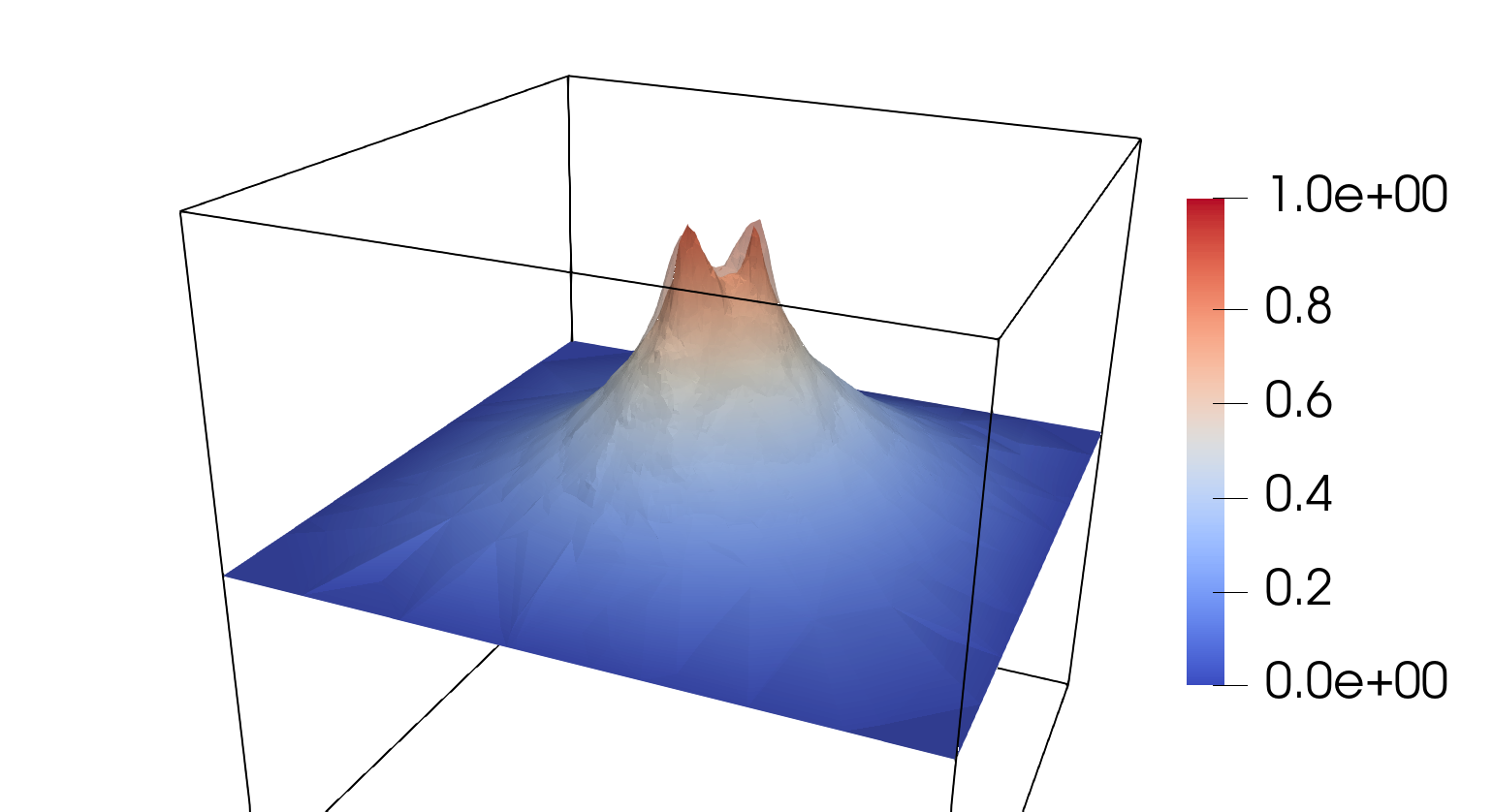}
		\caption{Reference solution (transparent, $N=43156$) and solution  for $\rho=0.1$ (opaque, $N=4574$).}
		\label{fig:warp-b}
	\end{subfigure}
	\caption{Test \ref{test_final}: Solutions on a plane orthogonal to the $z$-axis and located at $z=0$.}
	\label{fig:warp}
\end{figure}

Figure~\ref{fig:warp-a} shows the solutions obtained for $\rho=0$ and $\rho=0.1$ on a plane orthogonal to the $z$-axis and located at $z=0$. As expected, due to the non conformity of the chosen coarse mesh, standard finite element basis functions ($\rho=0$) are not able to capture the steep gradient close to the inclusion. Only by refining the mesh it is possible to reproduce the solution obtained for $\rho=0.1$. This is shown in Figure \ref{fig:warp-b}, where a good agreement between the solution obtained with $\rho=0.1$ and the reference solution can be observed.
We recall that the reference solution for this case is a solution computed with $\rho=0$ on a strongly refined mesh.

\begin{figure}
	\centering
	\begin{subfigure}{.45\textwidth}
		\centering
		\includegraphics[width=1.0\linewidth]{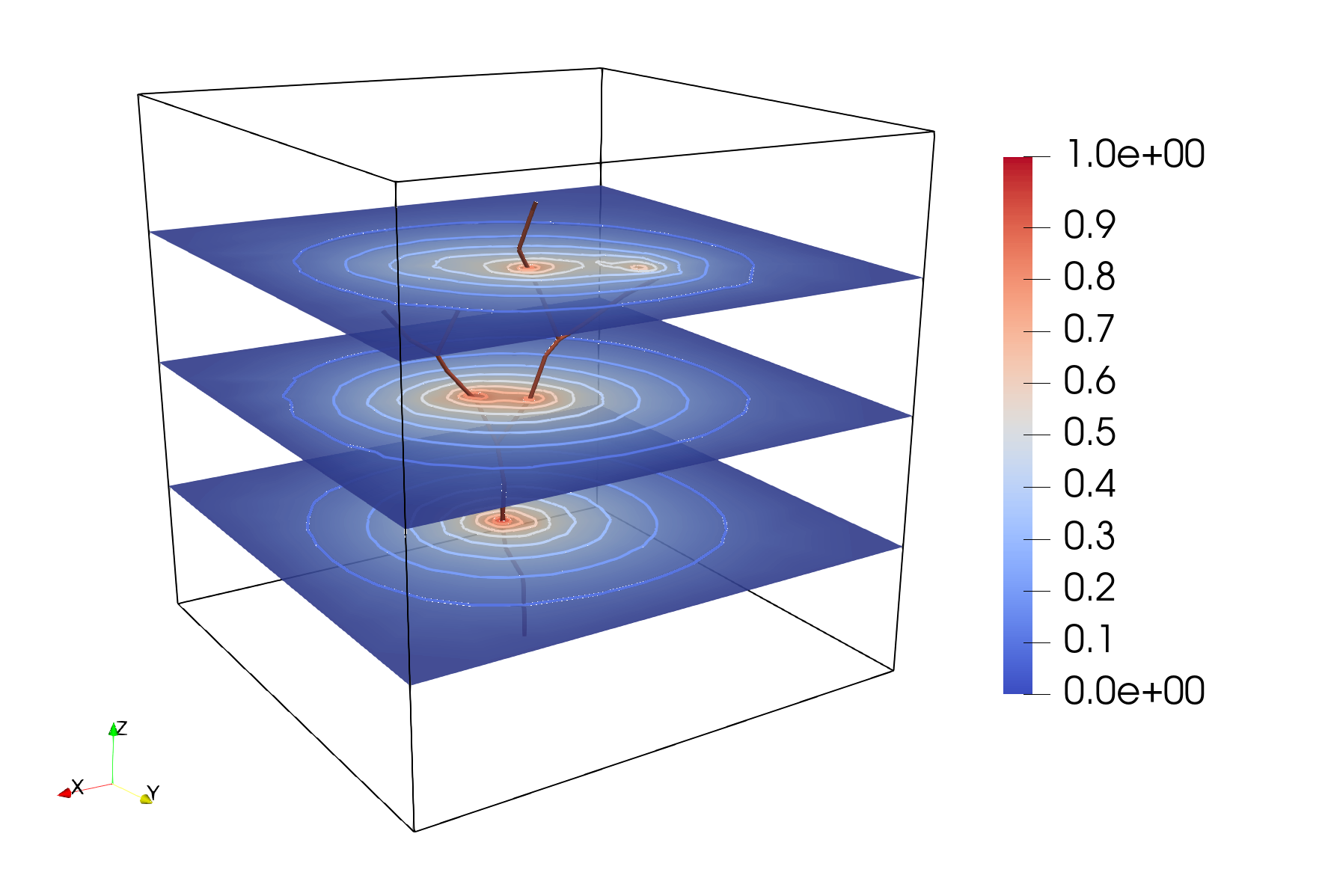}
		\caption{Solution for $\rho=0.1$ on sections orthogonal to the $z$-axis.}
		\label{fig:SolutionFEM_sliced-a}
	\end{subfigure}\hfill
	\begin{subfigure}{.45\textwidth}
		\centering
		\includegraphics[width=1.0\linewidth]{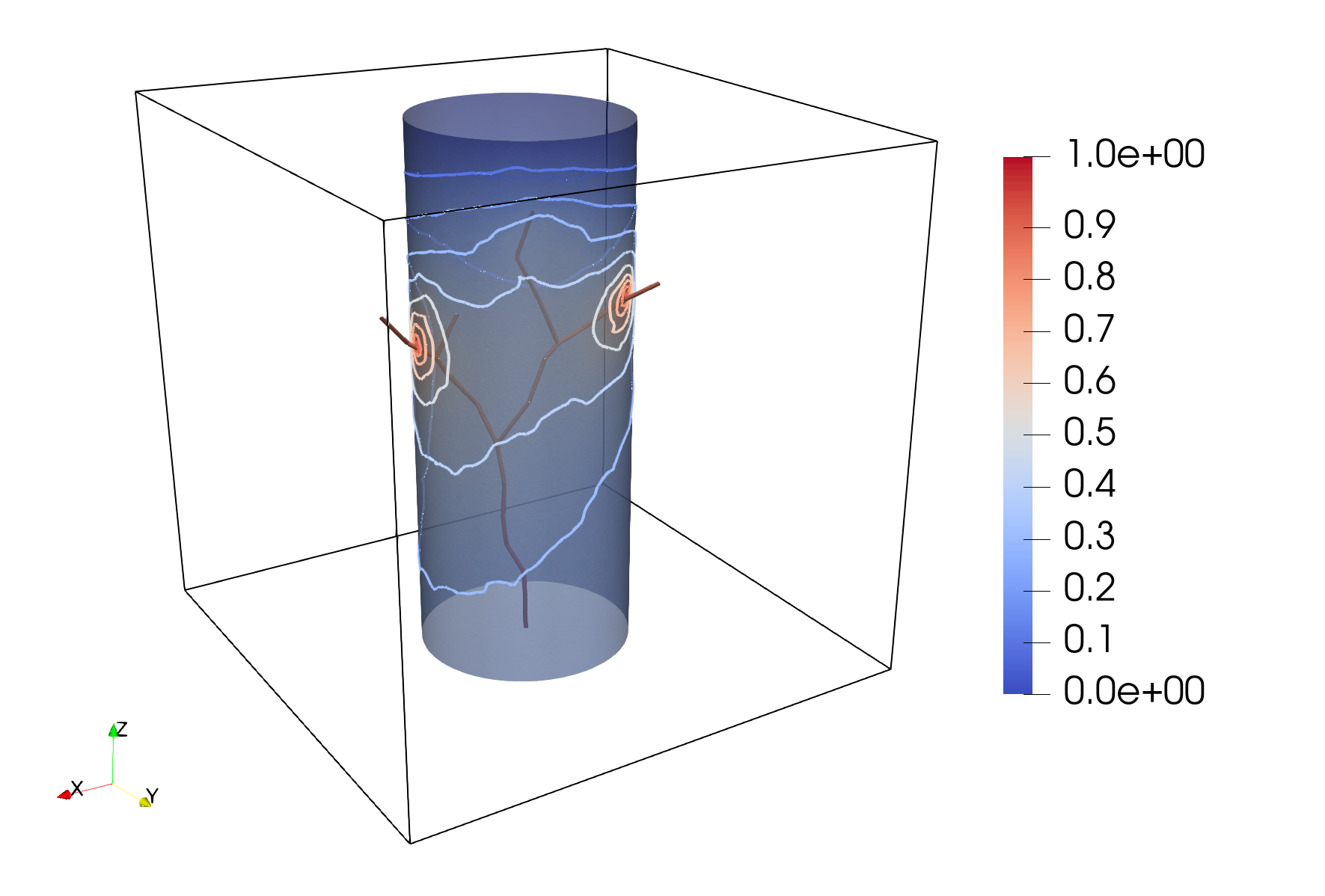}
		\caption{Solution for $\rho=0.1$ on a cylindrical surface of radius 0.4 and parallel to the $z$-axis. }
		\label{fig:SolutionFEM_sliced-b}
	\end{subfigure}
	\caption{Test \ref{test_final}: Solution obtained for $\rho=0.1$ on three different sections orthogonal to the $z$-axis and on a cylindrical surface parallel to the $z$-axis.}
	\label{fig:SolutionFEM_sliced}
\end{figure}

Finally, Figure~\ref{fig:SolutionFEM_sliced} reports the solution obtained for $\rho=0.1$ on three different sections of $\Omega$ orthogonal to the $z$-axis (Figure~\ref{fig:SolutionFEM_sliced-a}) and on a cylindrical surface parallel to the $z$-axis itself (Figure~\ref{fig:SolutionFEM_sliced-b}).

\section{Conclusions}
The present work presented an XFEM based implementation of a PDE-constrained optimization method for 3D-1D coupled problems. Suitable enrichment functions have been proposed to tackle problems with thin inclusions in very general cases: inclusion can have arbitrary orientations, form intersections or end inside the domain. A suitable quadrature rule has been introduced to numerically integrate the irregular enrichment functions on general polyhedral cells. The quadrature strategy uses the property of the enrichment functions of being regular in a direction tangential to inclusion centreline, and irregular, with a discontinuous gradient, on planes orthogonal to the centreline. Several numerical tests are proposed to show the effectiveness of the method in capturing the expected behavior of the solution also on meshes characterized by a maximum element diameter much larger than the radius of the inclusion. A validation of the methodology is also performed through the comparison with an available analytical solution, or with solutions obtained on adapted meshes.

\section*{Acknowledgements}
This publication is part of the project NODES which has received funding from the MUR-M4C2 1.5 of PNRR with grant agreement no. ECS00000036. Authors also acknowledge financial support from INdAM-GNCS.

\bibliographystyle{elsarticle-num}
\bibliography{xfem3D1D}
\end{document}